\documentclass[11pt,reqno]{amsart}

\usepackage{amsthm}
\usepackage{amssymb}
\usepackage{latexsym}
\usepackage{multicol}
\usepackage{verbatim,enumerate}
\usepackage{hyperref}
\usepackage{amsmath, amscd}
\usepackage{mathrsfs}
\usepackage[all,cmtip]{xy}
\usepackage{soul}

\advance\textwidth by 1.2in \advance\oddsidemargin by -.6in \advance\evensidemargin by -.6in
\usepackage{tikz}
\usetikzlibrary{decorations.markings}
\tikzstyle{vertex}=[ circle, fill, draw, inner sep=0pt, minimum size=4pt,]
\tikzstyle{edge}= [thick]

\newtheorem*{cor}{Corollary}%[section]
\newtheorem*{lem}{Lemma}
\newtheorem*{prop}{Proposition}

\theoremstyle{definition} \newtheorem*{defn}{Definition}
\theoremstyle{definition}
\newtheorem{thm}{Theorem}
\newtheorem*{thm*}{Theorem}

\newtheorem*{rem}{Remark}

\newenvironment{pf}{\proof}{\endproof}
\newcounter{cnt}
\newenvironment{enumerit}{\begin{list}{{\hfill\rm(\roman{cnt})\hfill}}{%
\settowidth{\labelwidth}{{\rm(iv)}}\leftmargin=\labelwidth%
\advance\leftmargin by \labelsep\rightmargin=0pt\usecounter{cnt}}}{\end{list}} \makeatletter
\def\mydggeometry{\makeatletter\dg@YGRID=1\dg@XGRID=20\unitlength=0.003pt\makeatother}
\makeatother \theoremstyle{remark}

\numberwithin{equation}{section}

 \DeclareMathOperator{\Ht}{ht}

\newcommand{\id}{\operatorname{id}}

\newcommand{\wt}{\operatorname{wt}}

\newcommand{\nc}{\newcommand}
\newcommand{\rnc}{\renewcommand}

\nc{\cal}{\mathcal} \nc{\goth}{\mathfrak} \rnc{\bold}{\mathbf}

\renewcommand{\Bbb}{\mathbb}
\nc\bomega{{\mbox{\boldmath $\omega$}}} \nc\bpsi{{\mbox{\boldmath $\Psi$}}}
 \nc\balpha{{\mbox{\boldmath $\alpha$}}}
 \nc\bbeta{{\mbox{\boldmath $\beta$}}}
 \nc\bpi{{\mbox{\boldmath $\pi$}}}
  \nc\bpis{{\mbox{\boldmath \scriptsize$\pi$}}}
 \nc\bullets{{\mbox{\scriptsize $\bullet$}}}
 \nc\bvarpis{{\mbox{\boldmath \scriptsize$\varpi$}}}
  \nc\bvarpi{{\mbox{\boldmath $\varpi$}}}

\nc\bepsilon{{\mbox{\boldmath $\epsilon$}}}
  \nc\bomegas{{\mbox{\boldmath\scriptsize $\omega$}}}
  \nc\bepsilons{{\mbox{\boldmath \scriptsize$\epsilon$}}}

\nc\hlien{\hat{\lie n}^+}
  \nc\bxi{{\mbox{\boldmath $\xi$}}}
\nc\bmu{{\mbox{\boldmath $\mu$}}} \nc\bcN{{\mbox{\boldmath $\cal{N}$}}} \nc\bcm{{\mbox{\boldmath $\cal{M}$}}} \nc\blambda{{\mbox{\boldmath
$\lambda$}}}%\nc\mathbb Nu{{\mbox{\boldmath $\nu$}}}

\newcommand{\lie}[1]{\mathfrak{#1}}

\makeatletter
\def\section{\def\@secnumfont{\mdseries}\@startsection{section}{1}%
  \z@{.7\linespacing\@plus\linespacing}{.5\linespacing}%
  {\normalfont\scshape\centering}}
\def\subsection{\def\@secnumfont{\bfseries}\@startsection{subsection}{2}%
  {\parindent}{.5\linespacing\@plus.7\linespacing}{-.5em}%
  {\normalfont\bfseries}}
\makeatother

 \nc{\Hom}{\operatorname{Hom}}
  \nc{\mode}{\operatorname{mod}}
\nc{\End}{\operatorname{End}} \nc{\wh}[1]{\widehat{#1}} \nc{\Ext}{\operatorname{Ext}}
 \nc{\ch}{\operatorname{ch}} \nc{\ev}{\operatorname{ev}}
\nc{\Ob}{\operatorname{Ob}} \nc{\soc}{\operatorname{soc}} \nc{\rad}{\operatorname{rad}} \nc{\head}{\operatorname{head}}
\def\Im{\operatorname{Im}}

 \nc{\Cal}{\cal} \nc{\Xp}[1]{X^+(#1)} \nc{\Xm}[1]{X^-(#1)}
\nc{\on}{\operatorname} \nc{\Z}{{\bold Z}} \nc{\J}{{\cal J}} \nc{\C}{{\bold C}} \nc{\Q}{{\bold Q}}

\nc{\N}{{\Bbb N}} \nc\boa{\bold a} \nc\bob{\bold b} \nc\boc{\bold c} \nc\bod{\bold d} \nc\boe{\bold e} \nc\bof{\bold f} \nc\bog{\bold g}
\nc\boh{\bold h} \nc\boi{\bold i} \nc\boj{\bold j} \nc\bok{\bold k} \nc\bol{\bold l} \nc\bom{\bold m} \nc\bon{\bold n} \nc\boo{\bold o}
\nc\bop{\bold p} \nc\boq{\bold q} \nc\bor{\bold r} \nc\bos{\bold s} \nc\boT{\bold t} \nc\boF{\bold F} \nc\bou{\bold u} \nc\bov{\bold v}
\nc\bow{\bold w} \nc\boz{\bold z} \nc\boy{\bold y} \nc\ba{\bold A} \nc\bb{\bold B} \nc\bc{\mathbb C} \nc\bd{\bold D} \nc\be{\bold E} \nc\bg{\bold
G} \nc\bh{\bold H} \nc\bi{\bold I} \nc\bj{\bold J} \nc\bk{\bold K} \nc\bl{\bold L} \nc\bm{\bold M}  \nc\bo{\bold O} \nc\bp{\bold
P} \nc\bq{\bold Q} \nc\br{\bold R} \nc\bs{\bold S} \nc\bt{\bold T} \nc\bu{\bold U} \nc\bv{\bold V} \nc\bw{\bold W} \nc\bx{\bold
x} \nc\KR{\bold{KR}} \nc\rk{\bold{rk}} \nc\het{\text{ht }}
\nc\bz{\mathbb Z}
\nc\bn{\mathbb N}

\nc\toa{\tilde a} \nc\tob{\tilde b} \nc\toc{\tilde c} \nc\tod{\tilde d} \nc\toe{\tilde e} \nc\tof{\tilde f} \nc\tog{\tilde g} \nc\toh{\tilde h}
\nc\toi{\tilde i} \nc\toj{\tilde j} \nc\tok{\tilde k} \nc\tol{\tilde l} \nc\tom{\tilde m} \nc\ton{\tilde n} \nc\too{\tilde o} \nc\toq{\tilde q}
\nc\tor{\tilde r} \nc\tos{\tilde s} \nc\toT{\tilde t} \nc\tou{\tilde u} \nc\tov{\tilde v} \nc\tow{\tilde w} \nc\toz{\tilde z} \nc\woi{w_{\omega_i}}
\nc\chara{\operatorname{Char}}
\begin{document}
\title[]{Higher order Kirillov--Reshetikhin modules for $\bu_q(A_n^{(1)})$,  Imaginary modules  and Monoidal Categorification  }
\author{Matheus Brito}
\address{Departamento de Matematica, UFPR, Curitiba - PR - Brazil, 81530-015}
\email{mbrito@ufpr.br}
\thanks{M.B. thanks the Hausdorff Research Institute for Mathematics and the
organizers of the Trimester Program “New Trends in Representation Theory” for excellent
working conditions.}
\author{Vyjayanthi Chari}
\address{Department of Mathematics, University of California, Riverside, 900 University Ave., Riverside, CA 92521, USA}
\email{chari@math.ucr.edu}
\thanks{V.C. was partially supported by DMS-1719357, the Max Planck Institute, Bonn and by the  InfoSys Visiting Chair position at the Indian Institute of Science.}
\thanks{Both authors were supported by the Oberwolfach Research Fellows program and acknowledge the superb working conditions at the Mathematisches Forschungsinstitut Oberwolfach.}

\begin{abstract} We study the  family of irreducible modules for quantum affine $\lie{sl}_{n+1}$ whose Drinfeld polynomials are supported on just one node of the Dynkin diagram.  We identify all  the prime modules in this family and prove a unique factorization theorem.  The  Drinfeld polynomials of the prime modules encode information coming from the points of reducibility of tensor products of the fundamental modules  associated to $A_m$ with $m\le n$. These prime modules are  a special class of the snake modules studied by Mukhin and Young.  We relate our modules  to the work of Hernandez and Leclerc and define  generalizations of the category $\mathscr C^-$. This leads naturally to the notion of an inflation of the corresponding Grothendieck ring.  In  the last section we show that the tensor product of a (higher order)  Kirillov--Reshetikhin module with its dual always contains an imaginary module in its Jordan--Holder series and give an explicit formula for its Drinfeld polynomial. Together with the results of \cite{HL13a} this gives examples of a product of cluster variables which are not in the span of cluster monomials. We also discuss the connection of our work with the examples arising from the work of \cite{LM18}.
  Finally, we  use our methods to  give a family  of imaginary modules in type $D_4$ which do not arise from an embedding of $A_r$ with $r\le 3$ in $D_4$. 
\end{abstract}

\maketitle
\section*{Introduction} In \cite{KR87}, Kirillov and Reshetikhin introduced a family of irreducible  finite--dimensional modules for the Yangian of a simple Lie algebra and conjectured a character formula for these representations. Subsequently, these modules were also defined for the quantum loop algebra associated to a simple Lie algebra and these are now called the Kirillov--Reshetikhin (KR) modules. These modules have many nice properties and there is an extensive literature on the subject, which includes their connections with integrable systems, the  combinatorics of   crystal bases, the fermionic formula \cite{FK08,FK13, HKOTY99, Ke11, Nak02,  OSS18, SW10}, and more recently \cite{HL10, HL13a, HL13} they have been shown to be   connected to   cluster algebras via the notion of monoidal categorification.\\\\
In this paper we are interested in higher order versions of KR--modules for a quantum loop algebra associated to a simple Lie algebra of type $A_n$.  To explain this, we recall that the KR--modules are indexed by three parameters $(i,a,r)$ where $i\in[1,n]$, $a$ is a complex number and $r$ a positive integer and the corresponding module is denoted $W_{i,a}^{(r)}$. Equivalently, one can think of these modules as being indexed by $i\in[1,n]$ and the $q$--segment $(a,aq^2,\cdots, aq^{2r-2})$. The appearance of $q^2$ in the segment can be explained as follows.
In the case $n=1$ it is known that the tensor product 
$W_{1,a}^{(1)}\otimes W_{1,b}^{(1)}$ is reducible if and only if $b=aq^{\pm 2}$ and also that the  module $W_{1,a}^{(r)}$ is the simple socle of the tensor product $W_{1,a}^{(1)}\otimes W_{1,aq^2}^{(1)}\otimes\cdots\otimes W_{1,aq^{2r-2}}^{(1)}$.
In the higher rank case, it remains true that $W_{i,a}^{(1)}\otimes W_{i,b}^{(1)}$ is reducible if $b=aq^{\pm 2}$ but there are many more points of reducibility. For instance in the case of $A_3$ the module $W_{2,a}^{(1)}\otimes W_{2,aq^{\pm 4}}^{(1)}$ is also reducible. A complete list of the points of reducibility simple Lie algebras of classical type can be found in \cite{Ch01} (see also \cite{AK97}). \\\\ 
By a higher order KR--module  for the quantum loop algebra of type $A_n$  we mean a module which is indexed by a pair $(i,\boa)$ where $\boa$ is an increasing $r$--tuple of integers $(a_1,\cdots, a_r)$, $r\ge 2$    such that  $W_{i,q^{a_j}}^{(1)}\otimes W_{i,q^{a_{j-1}}}^{(1)}$ is reducible and at least one of $a_j-a_{j-1}\ne 2$. We call the family which includes both the  KR--modules and higher order KR--modules at node $i$  as KR--modules associated to  $(i,n)$--segments. In fact these modules are very special examples of the prime snake modules which appear in the work of \cite{MY12a,MY12}. The  results of those two papers are used extensively in our work.\\\\ 
We describe the results of this paper. Following \cite{CP91} we introduce the notion of two $(i,n)$--segments being in general or special position. We prove that the associated tensor product is irreducible if and only if the segments are in pairwise general position. (In particular, we recover a special case of the result of \cite{MY12}, namely that the modules associated to an $(i,n)$--segments are prime, i.e.,  cannot  be written as a tensor product of nontrivial representations). This allows us to prove a classification result and a unique factorization result. Namely, suppose that we have an irreducible module for quantum affine $A_n$  whose associated Drinfeld polynomial (see \cite{CP94}) is trivial at all nodes different from $i$. Then it is prime if and only if the Drinfeld polynomial defines an $(i,n)$--segment. Otherwise it can be written uniquely as a tensor product of modules associated with $(i,n)$--segments in general position.\\\\
In Section \ref{secinflation}, we study the relation between the higher order KR--modules and the work of Hernandez and Leclerc. In particular given a triple of integers $(\bar i, i, n)$ with $n+1=i(\bar i+1)$ we define a full subcategory $\mathscr C_{\bar i, n}$ of finite--dimensional representations of the quantum loop algebra of type $A_n$. In the case when $\bar i=n$ this is precisely the category $\mathscr C^-$ defined in \cite{HL13a}.
We show that $\mathscr C_{\bar i, n}$ is a monoidal tensor category whose Grothendieck ring is isomorphic to the Grothendieck ring $\mathscr C_{\bar i,\bar i}$ (equivalently the Grothendieck ring of $\mathscr C^-$ for the quantum loop algebra associated to $A_{\bar i} $). The  isomorphism is defined by requiring that the class of the  module $W_{j,a}^{(1)}$  for $A_{\bar i}$ maps to the class of  $W_{ij, ia}^{(1)}$ for $A_n$.  We conjecture that in general the isomorphism maps the class of an irreducible object of $\mathscr C_{\bar i, \bar i}$ to the class of an irreducible object of $\mathscr C_{\bar i, n}$. We prove that the conjecture holds for any KR--module associated to a $(j,\bar i)$--segment. We also establish the conjecture for the category $\mathscr C_1$ defined in \cite{HL10}. For this we define the analogous category $\mathscr C_{\bar i, n}^1$ and show that it is a tensor category and is the image of $\mathscr C_{\bar i, \bar i}^1.$ It now   follows from \cite{HL13} that  $\cal K_0(\mathscr C_{\bar i, n}^1)$ is a monoidal categorification of the cluster algebra of type $A_{\bar i}$ with the frozen variables being higher order KR--modules.\\\\
Recall that  a finite--dimensional irreducible module is said to be imaginary if its tensor square is reducible.
The first example of such a  module was given by Leclerc in \cite{Lec02} in types $A_4$, $B_3$, $C_2$ and $G_2$. Using embeddings of the associated quantum loop algebra it was immediate that they existed in all simple Lie algebras of higher rank. This was the only example in the literature of an imaginary modules  until  the work  \cite{LM18} of Lapid and Minguez which gives examples in type $A_n$. Their work, as well as the original example, come from affine Hecke algebras and affine Schur--Weyl duality. 
\\\\
In Section \ref{imaginary} of this  paper we use the results of the earlier sections to  give a systematic way to construct imaginary modules. We work with tensor products of dual KR--modules associated to $(i,n)$--segments. In particular, we prove (see Theorem \ref{imagen} for the  general statement) that if  $r\ge 2$ the tensor product  $W_{i,a}^{(r)}\otimes W_{n+1-i,a-n-1}^{(r)}$   has an imaginary module in its Jordan--Holder series. The example  given by  Leclerc in type $A_4$   can be viewed as coming from our result  in type $A_3$ and the tensor product $W_{2,a}^{(2)}\otimes W_{2,aq^4}^{(2)}$. Our examples do not satisfy the conditions for imaginary modules imposed in  the work of  \cite{LM18} except  for small values of $r$.  In type $D_4$ we find an imaginary module in the tensor product of $W_{1,a}^{(r)}$ and its dual for all $r\ge 3$. This example is not obtained by an embedding of any $A_r$--subalgebra in $D_4$, in fact the associated Drinfeld polynomial  would define a real module for any $A_3$ type subalgebra.  \\\\
We conclude the introduction with a brief discussion on higher order KR--modules in other types. They clearly exist, for instance in type $D_4$ one could just look at the top irreducible constituent  of the tensor product of a $W_{i,a}^{(1)} $  and its dual. This representation is obviously prime. However the problem arises in these higher ranks because the points of reducibility of $W_{i,a}^{(1)}\otimes W_{i,b}^{(1)}$ have \lq\lq gaps\rq\rq. In $D_n$ and $i=1$ the reducibility is when $a=bq^{\pm 2}$ and $a=bq^{2n-2}$. (In type $A_n$ in contrast the points of reducibility are of the form $\{aq^{\pm 2s}: 1\le s\le p\}$ for a suitable choice of $p$). The existence of the gaps makes it difficult to define a good analog of $(i,n)$--segments; it would appear that one has to \lq\lq glue\rq\rq \; two different kinds of segments together multiple times. The difficulty is analogous in some sense to the difficulty in type $A_2$ of gluing segments associated to nodes 1 and 2 to give a prime representation. We hope to return to these ideas in the future.
\\\\
{\em Acknowledgements. The authors thank David Hernandez and Bernard Leclerc for many interesting discussions on their work over the years. The authors also thank the referee for their careful reading of the paper and the many helpful comments and remarks. In particular, the referee has pointed out that there might be analogous statements for the generalizations in \cite{FHOO} of the category $\mathscr C^-$ and we  hope to return to this in future work. }
\section{Preliminaries}
In this section we   set up the notation to be used in the paper and recall several results which will play a crucial role in the later sections.
\\\\
We assume throughout that $q$ is a non--zero complex number and not a root of unity. As usual $\mathbb C$ (resp. $\mathbb C^\times$, $\mathbb Z$, $\mathbb Z_+$, $\mathbb N$) will denote the set of complex numbers (resp. non-zero complex numbers, integers, non-negative integers, positive integers). Given $i, j\in \mathbb N$ with $i\le j$ we let $[i,j]$ be the set of integers $\{i, i+1,\cdots, j\}$. 

\subsection{The algebra $\hat\bu_n$}\label{basicdef}  Let $\hat\bu_n$  be  the quantum loop algebra over $\mathbb C$ associated to $\lie{sl}_{n+1}$;  we refer the reader to \cite{CP94} for precise definitions. For our purposes, it is enough to recall that $\hat\bu_n$ is a Hopf algebra and is generated as an algebra  by elements $x_{i,s}^\pm$, $\phi^\pm_{i,s}$ $i\in[1,n]$ and $s\in\mathbb Z$. The algebra  generated by the elements $\phi^\pm_{i,s}$, $i\in[1,n]$,  $s\in\mathbb Z$, is denoted by $\hat\bu_n^0$ and is a   commutative subalgebra of $\hat\bu_n$.\\\\ Given  $a\in\mathbb C^\times$ let $\tau_a:\hat\bu_n\to\hat\bu_n$ be the Hopf algebra homomorphism given by $$x_{i,s}^\pm\to a^s x_{i,s}^\pm, \ \ \tau_a(\phi_{i,s}^\pm)=a^s\phi^\pm_{i,s}, \ \ i\in[1,n],  \ s\in\mathbb Z.$$ 
The quantum affine  analog of the Cartan involution of $A_n$  is the unique algebra involution $\Omega:\hat\bu_n\to \hat\bu_n$ given by $$\Omega(x^\pm_{i,s})=-x^\mp_{i,-s},  \ \ \Omega(\phi_{i,s}^\pm ) = \phi_{i,-s}^\mp,\ \ i\in[1,n],\  s\in\mathbb Z.$$ \\\\
Given  $J=[i,j]\subset [1,n]$ let $\hat{\bu}_{n,J}$ be the subalgebra of $\hat\bu_n$ generated by the elements $x_{p,s}^\pm, \phi_{p,s}^\pm$, $s\in\mathbb Z$,  $p\in J$. Then we have an  algebra (but not a Hopf algebra) homomorphism  $\hat\bu_{j-i+1}\to\hat\bu_{n,J}$.
\subsection{The $\ell$--weight and $\ell$--root lattice} Let $\cal P_n$ (resp. $\cal P_n^+$) be the (multiplicative)  free abelian group (resp. monoid) generated by  elements  $\{\bomega_{i,a}: i\in [1,n],\ \ a\in\mathbb Z\}$  and  denote by $\bold 1$ the identity element.
It will be convenient to set $\bomega_{k,a}=\bold 1$ for all $a\in\mathbb Z$ if $k\notin[1,n]$. The elements of $\cal P_n$ are called $\ell$--weights and those of $\cal P_n^+$ the dominant $\ell$--weights.  Let $P_n$ be the free (additive) abelian group on generators $\{\omega_i:i \in [1,n]\}$ and $P_n^+$ the corresponding monoid. Define a morphism of groups by extending the assignment $$\wt:\cal P_n\to P_n,\ \ \wt\bomega_{i,a}=\omega_i, \ \ 
i\in [1,n], \ a\in\mathbb Z.$$ 
For $i\in[1,n]$ and $a\in\mathbb Z$ define  $\balpha_{i,a}\in\cal P_{n}$ and $\alpha_i\in P_n$ by
$$\balpha_{i,a}=\bomega_{i-1,a}^{-1}\bomega_{i, a-1}\bomega_{i,a+1}\bomega_{i+1,a}^{-1},\ \ \alpha_i=2\omega_i-\omega_{i-1}-\omega_{i+1}$$ and let $\cal Q_{n}$ (resp. $\cal Q_n^+$) be the subgroup (resp. submonoid) of $\cal P_n$ generated by these elements. The subgroup (resp. submonoid)  of $P_n$ spanned by the elements $\alpha_i, i\in[1,n]$, is denoted $Q_n$ (resp. $Q_n^+$).\\\\
Given  $J=[i,j]\subset[1,n]$ let $\cal P_{n,J}$ be the subgroup of $\cal P_n$ generated by the elements $\bomega_{j,c}$ with $j\in J$ and $c\in\mathbb Z$. Clearly $\cal P_{n,J}\cong \cal P_{j-i+1}$ and we shall use this freely without mention. Define $\cal P_{n,J}^+$  in the obvious way and identify it with the corresponding monoid in $\hat\bu_{j-i+1}$.\\\\
Define a homomorphism  $\cal P_n\to\cal P_{j-i+1}$ sending $\bomega\mapsto\bomega_J$    by extending the assignment   $$\bomega_{i,c}\mapsto\bomega_{i,c},\ \ i\in J,\ \ \bomega_{i,c}\mapsto \bold 1,\ \ i\notin J.$$ Clearly the homomorphism maps $\cal Q_n$ to $\cal Q_{j-i+1}$.

\subsection{ The category $\mathscr F_n$}\label{fn}  Denote by   $\mathscr F_{n}$   the category of type 1 finite--dimensional representations of $\hat\bu_n$. In particular if $V$ is any object of $\mathscr F_n$ we can write $$V=\bigoplus_{\mu\in P} V_\mu,\ \ V_\mu=\{v\in V: \phi_{i,0}^{\pm}v=q^{\pm \mu_i}v\},\ \ \mu=\sum_{i=1}^n\mu_i\omega_i\,\ \ \mu_i\in\mathbb Z,\  i\in[1,n],$$ and we set $\wt V=\{\mu\in P: V_\mu\ne 0\}.$ Given any object $V$ of $\mathscr F_{n}$ we can regard it as a module for the commutative subalgebra $\hat{\bu}_{n}^0$. It follows that we can write $V$ as a direct sum of generalized  eigenspaces for the action of this subalgebra. The generalized eigenspaces are called $\ell$--weight spaces.\\\\
The Hopf algebra structure on $\hat\bu_n$ is 
chosen so that the comultiplication on the generators $x_{i,0}^+$  is given by $\Delta(x_{i,0}^+) = x_{i,0}^+\otimes 1 + \phi_{i,0}^+\otimes x_{i,0}^+$ (see \cite{BCKV} for details). 
Note that the Hopf algebra structure ensured that $\mathscr F_n$ contains the trivial representation, and is closed under tensor products and duals. \\\\
\begin{defn}
    We say that an object $V$ of $\mathscr F_n$ is prime if it is not isomorphic to $U\otimes U'$ where $U$ and $U'$ are non--trivial objects of $\mathscr F_n$.\hfill\qedsymbol
\end{defn}
\noindent Let $\cal K_0(\mathscr F_n)$ be  the corresponding  Grothendieck ring of $\mathscr F_n$ and denote by  $[V]$  the isomorphism class  of an object $V$ of $\mathscr F_n$.\\\\
Any object $V$ of $\mathscr F_n$  has a right and a  left dual  denoted by $V^*$ and ${}^*V$ respectively, and we have $\hat\bu_n$--maps
 $$\mathbb C\hookrightarrow V^*\otimes V,\ \ {}^*V\otimes V\to \mathbb C\to 0.$$
We shall freely use properties of duals, in particular, the isomorphisms
\begin{gather*}(W\otimes V)^*\cong V^*\otimes W^*,\ \  {}^*(W\otimes V)\cong {}^*V\otimes {}^*W,\\
 \Hom_{\hat \bu_n}(V\otimes U, W)\cong \Hom_{\hat \bu_n}( U, V^*\otimes W),\ \ \Hom_{\hat \bu_n}(U\otimes V, W)\cong \Hom_{\hat\bu_n}( U, W\otimes {}^*V).\end{gather*}
For $a\in\mathbb C^\times$ and $V\in\mathscr F_n$, let $\tau_a V$ and $\Omega(V)$ be the corresponding objects of $\mathscr F_n$ obtained by pulling $V$ back by the automorphisms $\tau_a$ and $\Omega$ respectively. Then  $$\tau_a(V\otimes W)\cong \tau_a V\otimes \tau_a W,\ \ \Omega(V\otimes W)\cong \Omega(W)\otimes \Omega(V).$$

 \subsubsection{The modules $W(\bomega)$ and $V(\bomega)$ }\label{weyldef} It is  convenient to  identify $ \cal P_n^+$ with the monoid consisting of   $n$--tuple of polynomials by extending the assignment  $\bomega_{i,a}\mapsto (1-\delta_{i,j}q^au)_{j\in [1,n]}$ to a multiplicative homomorphism. For  $\bomega\in\cal P_n^+$  we let $W(\bomega)$ be the $\hat\bu_n$--module generated by an element $v_\bomegas$ satisfying the relations
 $$x_{i,s}^+ v_{\bomegas} = 0 = (x_{i,0}^-)^{\deg\pi_i(u)+1}v_\bomegas, \ %\ k_iv_\bomegas = q^{\deg\pi_i(u)}v_\bomegas, \ 
 \    \phi_{i,s}^\pm v_{\bomegas} = \gamma_{i,s}^\pm v_\bomegas,\ \ 
 i\in [1,n], \ s\in\mathbb{Z},$$  where $\gamma_{i,s}^\pm \in \mathbb C$ are defined by  $$ \sum_{s=0}^\infty \gamma_{i,\pm s}^{\pm}u^{\pm s} = q^{\deg \pi_i}\frac{\pi_i(q^{-1} u)}{\pi_i(qu)}, \ \ \ \ \bomega = (\pi_i(u))_{i\in I}.$$ 
 Moreover $$\dim W(\bomega)_{\wt\bomegas}=1,\ \ {\rm and\ for} \  \mu\in P,\ \  \dim W(\bomega)_\mu\ne 0\implies \wt\bomega-\mu\in Q^+.$$
Any quotient of $W(\bomega)$ is called an $\ell$--highest weight module  with highest $\ell$--weight $\bomega$ and we  continue to denote  by $v_\bomegas$ the image of the generator of $W(\bomega)$ in any quotient. It follows from \cite{CP01} that the module $W(\bomega)$ is finite--dimensional and has a unique irreducible quotient which we denote as $V(\bomega)$. 
Finally, any irreducible module in $\mathscr F_n$ is isomorphic to  a tensor product of objects of the form  $\tau_b V(\bomega)$ for some $b\in\mathbb C^\times$ and  $\bomega\in\cal P^+_n$.  \\\\
Setting  $$\bomega_{i,a}^*=\bomega_{n+1-i, a+n+1}
,\ \ {}^*\bomega_{i,a}=\bomega_{n+1-i,a-n-1},\ \ \Omega(\bomega_{i,a})=\bomega_{n+1-i, -a} $$  we get  corresponding isomorphisms $\bomega\mapsto\bomega^*$, $\bomega\mapsto {}^*\bomega$  and $\bomega\mapsto\Omega(\bomega)$ of $\cal P_n$ and $$V(\bomega)^*\cong V(\bomega^*),\ \ {}^*V(\bomega)\cong V({}^*\bomega), \ \  \Omega(V(\bomega))\cong V(\Omega(\bomega)).$$
\subsubsection{Tensor products}
 Part (i)  of the next proposition was proved in \cite{CP94, CP95}. Part (ii) was proved independently in \cite{AK97} and \cite{Ch01}  and part (iii) in \cite{Ch01}. Given $m,r\in[1,n]$ set
  \begin{equation}\label{smkn} S_{m,r,n}=\{2p+2-m-r: \max\{m,r\}\le p\leq \min\{m+r-1,n\}\}.\end{equation}
 \begin{prop}\label{weyl} Suppose that $\bomega=\bomega_{i_1,a_1}\cdots\bomega_{i_k,a_k}\in\cal P^+_n$ with $a_1\le \cdots\le a_k$.
 \begin{enumerit}
 \item[(i)] Let $\bomega'\in\cal P^+_n$. The module $V(\bomega\bomega')$ occurs with multiplicity one in the Jordan--Holder series of $V(\bomega)\otimes V(\bomega')$. Moreover, $V(\bomega\bomega')$ is isomorphic to $V(\bomega)\otimes V(\bomega')$ if and only if $V(\bomega)\otimes V(\bomega')$ and its left (or right) dual are $\ell$--highest weight modules.
 \item[(ii)]  We have, $$W(\bomega)\cong V(\bomega_{i_k,a_k})\otimes\cdots\otimes V(\bomega_{i_1,a_1})$$ and hence for all $\bomega,\bomega'\in\cal P^+_n$ the following holds in $\cal K_0(\mathscr F_n)$: $$[W(\bomega\bomega')]=[W(\bomega)][W(\bomega')].$$
 \item[(iii)] The module $V(\bomega_{i_1,a_1})
 \otimes V(\bomega_{i_2,a_2})\otimes \cdots \otimes V(\bomega_{i_k, a_k})$  is $\ell$--highest weight if and only if $a_j-a_{j-1}\notin S_{i_j,i_{j-1},n},\  \ 2\le j\le k$. If $a_{j-1}-a_j$ is also not in this set for all $2\le j\le k$ then the module is irreducible.
 \hfill\qedsymbol
 \end{enumerit}
 \end{prop}

\subsubsection{$\ell$--lowest weight modules }\label{lowestwt}  An $\ell$--lowest weight module is defined in the obvious way; it is generated by an element $v$ which is an eigenvector for the elements $\phi_{i,s}^\pm$ and $x_{i,s}^-v=0$ for all $i\in[1,n]$, $s\in\mathbb Z$.
\begin{prop} 
\begin{enumerit}
\item[(i)] 
Any $\ell$--highest weight module with $\ell$--highest weight $\bomega$ in $\mathscr F_n$ is also a lowest $\ell$--weight module with lowest weight $(\bomega^*)^{-1}$.
\item[(ii)] Let $V, V'$ be $\ell$--highest weight modules with $\ell$--highest weight  $\bomega,\bomega'\in\cal P^+_n$, respectively. Let   $v^-$ and $v^+$ be non--zero  lowest and highest $\ell$--weights of $V$ and $V'$. Then  $v^-\otimes v^+$ is an $\ell$--weight vector with $\ell$--weight $(\bomega^*)^{-1}\bomega'$ and $$V\otimes V'=\hat\bu_n(v^-\otimes v^+).$$ In particular if $U$ is a proper quotient of $V\otimes V'$ then $\dim U_{(\bomega^*)^{-1}\bomega'}\ne 0$.
\end{enumerit}
\end{prop} 
\begin{pf}
We sketch a proof. Let $\lambda=\wt\bomega$. Since $V$ is an $\ell$--highest weight module we have $\wt V\subset\lambda- Q^+$ and $\dim V_\lambda=1$.  Since $V$ is a finite--dimensional module for $\hat\bu_n$ and hence also for $\bu_n$ (the subalgebra generated by the elements $x_{i,0}^\pm,\  \phi_{i,0}^{\pm},\  i\in[1,n]$)  it follows that $\dim V_{w_\circ\lambda}=1$ where $w_\circ$ is the longest element of the Weyl group $S_{n+1}$ of $A_n$; in particular any non--zero element of  $V_{w_\circ\lambda}$ is an $\ell$--weight vector. It was shown in \cite{Ch01} that if $V=V(\bomega)$ then $V_{w_\circ\lambda}$ was an $\ell$--weight space with $\ell$--weight $(\bomega^*)^{-1}$. Since $V(\bomega)$ is a quotient of any $\ell$--highest weight module with $\ell$--weight $\bomega$, part (i) follows.\\\
Part (ii) is immediate from the formulae for the comultiplication \cite{Da98} (see also  \cite{Ch01}).
\end{pf}

\subsubsection{} We shall use the following consequence of Proposition \ref{lowestwt}.
\begin{lem}\label{socle} Suppose that $\bpi,\bpi_1,\bpi_2\in\cal P_n^+$. Then$$\Hom_{\hat\bu_n}(W(\bpi), V(\bpi_1)\otimes V(\bpi_2))\ne 0\implies \bpi\bpi_1^{-1}\in\wt_\ell V(\bpi_2)\ \ {\rm{and}}\ \ (\bpi^*)^{-1}\bpi_2^*\in\wt_\ell V(\bpi_1).$$\end{lem}
\begin{pf} Using duals we see that
$\Hom_{\hat\bu_n}(W(\bpi), V(\bpi_1)\otimes V(\bpi_2))\ne 0$ implies 
$$\Hom_{\hat\bu_n}(V(^*\bpi_1)\otimes W(\bpi),  V(\bpi_2))\ne 0,\ \ \Hom_{\hat\bu_n}(W(\bpi)\otimes V(\bpi_2)^*, V(\bpi_1)\ne 0.$$ Since $(^*\bpi_1)^*=\bpi_1$ it follows from Proposition \ref{lowestwt} that $$\bpi_1^{-1}\bpi\in\wt_\ell V(\bpi_2),\ \ (\bpi^*)^{-1}\bpi_2^*\in\wt_\ell V(\bpi_1).$$ 
\end{pf}
\subsection{The $\ell$--weight space decomposition and  $q$--characters}\label{qchar}  Let $\mathscr F_{n,\mathbb Z}$ be the full subcategory of $\mathscr F_n$ whose Jordan--Holder constituents are of the form $V(\bomega)$, $\bomega\in\cal P_n^+$. It is well--known that $\mathscr F_{n,\mathbb Z}$ is a rigid tensor subcategory of $\mathscr F_n$ and we let $\cal K_0(\mathscr F_{n,\mathbb Z})$ be the corresponding Grothendieck ring.
   It was proved in  \cite{FR99} that if $V$ is an object of $\mathscr F_{n,\mathbb Z
}$ then one can write $V$ as a direct sum of generalized eigenspaces for $\hat{\bu}_n^0$ and that the eigenvalues are indexed by elements of $\cal P_n$, $$V=\bigoplus_{\bomegas\in\cal P_n} V_\bomegas,\ \ \wt_\ell V=\{\bomega\in\cal P_n: V_\bomegas\ne 0\},\ \ \wt_\ell ^+ V=\wt_\ell V\cap \cal P^+_n.$$ Given a subgroup $\cal G$ of $\cal P_n$ we define $\chi^{\cal G}(V)$ to be the element of the group ring of $\cal G$ given by, \begin{equation}\label{truncg} \chi^{\cal G}(V)=\sum_{\bomegas\in \cal G}\dim V_\bomegas e(\bomega)\in\mathbb Z[\cal G].\end{equation}  The $q$--character of $V$ is the  element $\chi^{\cal P_n}(V)$.
\subsubsection{} The following was proved in \cite{FR99}. 
\begin{thm*} \label{fr}
The assignment $[V]\mapsto\chi^{\cal P_n}(V)$ gives an injective homomorphism   $\chi^{\cal P_n}: \cal K_0(\mathscr F_{n,\mathbb Z})\to \mathbb Z[\cal P_n]$ of rings. Moreover, $\cal K_0(\mathscr F_{n,\mathbb Z})$ is a polynomial ring in the  generators $[V(\bomega_{i,a})]$ with $i\in [1,n]$ and $ a\in\mathbb Z$. 
In particular if $V$ and $V'$ are objects of $\mathscr F_{n,\mathbb Z}$ we have, $$\wt_\ell (V\otimes V')=\wt_\ell V\wt_\ell V'. $$
\hfill\qedsymbol\end{thm*}
\subsubsection{} We shall use the following result which proof can be found in \cite{CM04}, for instance.
\begin{prop} Let $V$ be an $\ell$--highest weight module with $\ell$--highest $\bomega$. Then $$\wt_\ell V\subset\bomega (\cal Q^+)^{-1}.$$\hfill\qedsymbol  

\end{prop}

\subsection{Restrictions to $\hat\bu_{n,J}$}\label{diagsub}

  Given $J=[i,j]\subset [1,n]$ and $\bomega\in\cal P^+_n$ we have an isomorphism of $\hat{\bu}_{n,J}$--modules $V(\bomega_J)\cong\hat\bu_{n,J}v_\bomegas\subset V(\bomega)$ and $$V(\bomega_J)\otimes V(\bomega'_J)\cong \hat\bu_{n,J}v_\bomegas\otimes \hat\bu_{n,J}v_{\bomegas'}.$$ In particular if $v\in V(\bomega_J)\otimes V(\bomega'_J)$ is an $\ell$--highest weight vector with $\ell$--highest $\bomega_J\bomega_J'\balpha_J^{-1}$ where $\balpha\in\cal Q^+_{n,J}$, then $V(\bomega)\otimes V(\bomega')$ has  an $\ell$--highest weight vector with $\ell$--highest $\bomega\bomega'\balpha^{-1}$.

\subsection{Some results of Mukhin and Young}\label{MYsec}
We  recall some  results of Mukhin and Young which were established in \cite{MY12a, MY12} and which will play an important role in the subsequent sections.

\subsubsection{The set $\mathbb P_{i,a}$} \label{defpath}
For $i\in [1,n]$ and $a\in\mathbb Z$,  let $\mathbb P_{i,a}$ be the set of all functions $p: [0,n+1]\to\mathbb Z$ satisfying the following:
\begin{gather*} p(0)=i+a,\ \ p(r+1)-p(r)\in\{-1,1\},\ \ 0\le r\le n,\ \ p(n+1)=n+1-i+a.\\ \end{gather*}
For $p\in\mathbb P_{i,a}$ set  \begin{gather*}\boc_p^\pm =\{r\in[1,n]: p(r-1)=p(r)\pm 1=p(r+1)\},\\\\
\bomega(p)=\prod_{r\in \boc_p^+}\bomega_{r, p(r)}\prod_{r\in \boc_p^-}\bomega_{r,{p(r)}}^{-1}\in\cal P_n.\end{gather*}
In particular $\bomega(p)$ is in the subgroup  of $\cal P_n$ generated by the elements $\{\bomega_{j,c}: j\in[1,n],\ a\le c\le a+n+1 \}$.\\\\
Let $p_{i,a}$ and $p_{i,a}^*$ be the elements of $\mathbb P_{i,a}$ given as follows:
$$p_{i,a}(j)= \begin{cases} i-j+a,\ \ 0\le j\le i, \\ j-i+a,\ \ i<j\le n+1,\ \ 
\end{cases}\  p_{i,a}^*(j)=\begin{cases} a+i+j,\ \ 0\le j\le n+1-i,\\
a+2n+2-i-j,\ \ n+2-i\le j\le n+1.\end{cases}$$ Then $$\bomega(p_{i,a})=\bomega_{i,a},\ \ \bomega(p_{i,a}^*)=\bomega_{n+1-i, a+n+1}^{-1}.$$ 
The following is a simple calculation.
\begin{lem}\label{pims}
Let $a,b,c$ be integers  with  $b-a=2m_1$ and $c-b=2m_2$ for some $m_1, m_2\in\mathbb N$. Then, for all $p\in\mathbb P_{i,b}$ and  $j\in[0,n+1]$ we have  $$  p_{i,a}(j)<p(j)<p_{i,c}^*(j), \ \  p\in\mathbb P_{i,b}.$$\hfill\qedsymbol\end{lem}
\subsubsection{} \label{gjm}  More generally, given $j\in[1,n]$,  $m\in[0,\min\{j, n+1-j\}]$  and $a\in\mathbb Z$ it is not hard to see that there exist elements $g_{j,a}^m\in\mathbb P_{j,a}$ satisfying, $$\bomega(g_{j,a}^m)=\begin{cases} \bomega_{j,a},\ \ m =0,\\
\bomega_{n+1-j-m, a-n-1+2j+m} \bomega_{n+1-j, a-n-1+2j+2m}^{-1}\bomega_{j+m, a+m}, \  \ \ 2j>n+1,\\
\bomega_{j-m, a+m} \bomega_{n+1-j,a+n+1-2j+2m}^{-1}\bomega_{n+1-j+m,a+n+1-2j+m},\ \ \ \ 2j\le n+1.
\end{cases} $$ Notice that $$m=\min\{j,n+1-j\}\implies \bomega(g_{j,a}^m)=\bomega_{n+1-j, a+n+1}^{-1}=p_{j,a}^*.$$ It is easily checked that \begin{equation}\label{singleton}
\{g_{j,a}^m:1\le m\le \min\{j, n+1-j\}\}=\{g\in\mathbb P_{j,a}: \boc_g^-=\{n+1-j\}\}.\end{equation}
Similarly there exist elements $p_{j,a}^m\in\mathbb P_{j,a}$ such that 
$$\bomega(p_{j,a}^m)=\begin{cases}\bomega_{j,a},\ \ m=0,\\
\bomega_{j-m, a+m}\bomega_{j,a+2m}^{-1} \bomega_{j+m, a+m},\ \  0<m\leq \min\{j,n+1-j\}. \end{cases} $$and\begin{equation}\label{pjm}\{p_{j,a}^m: 1\le m\le 
\min\{j,n+1-j\}\}=\{p\in\mathbb P_{j,a}: \boc_p^-=\{j\}\}.\end{equation}
 \subsubsection{The set $\mathbb P_{\bomegas}$}\label{snakepaths}
Given  $\bomega=\bomega_{i_1,a_1}\cdots\bomega_{i_k,a_k}\in\cal P_n^+$ with $ a_1\le a_2\le \cdots\le a_k,$ define  $$\mathbb P_\bomegas\subset  \mathbb P_{i_1,a_1}\times\cdots\times \mathbb P_{i_k,a_k}$$ to consist of $k$--tuples  $(p_1,\cdots, p_k)$ satisfying: 
\begin{equation}\label{pbomega}  p_j(k)<p_s(k)\ {\rm{for \ all}}\  k\in[0,n+1]\ \ {\rm{and\ all}}\ \  1\le j<s\le r.\end{equation}
Given $\underline p\in\mathbb P_\bomegas$, set
\begin{equation}\label{omp}\bomega(\underline p)=\bomega(p_1)\cdots\bomega(p_r),\ \ \underline p=(p_1,\cdots, p_r).\end{equation}
The restriction in \eqref{pbomega} guarantees that  the expression on the right hand side of  \eqref{omp} is a reduced word in $\cal P_n$. Here we emphasize that this is equivalent to saying that there are no cancellations between the $\bomega(p_j)$. As an example, consider $n=3$ and $\bomega=\bomega_{2,0}\bomega_{2,4}$. Then $(p_{2,0}^*, p_{2,4})\notin \mathbb P_{\bomegas}$, since $p_{2,0}^*(2) = 4= p_{2,4}(2)$.  Moreover $$\bomega(p_{2,0}^*)\bomega(p_{2,4})= \bomega_{2,4}^{-1}\bomega_{2,4}=\bold 1\notin \wt_\ell V(\bomega).$$  
\\\\
For $i\in[1,n]$ and integers  $a_1\le \cdots\le a_k$, we shall write \begin{equation}\label{omegaia}\bomega_{i,\boa}=\bomega_{i,a_1}\cdots\bomega_{i,a_r},\ \  \mathbb P_{i,\boa}=\mathbb P_{\bomegas_{i,\boa}}.\end{equation}

\noindent Following \cite{MY12} we shall say that $\bomega\in\cal P_n^+$ is a prime snake  if we can write  \begin{equation}\label{ps} \bomega=\bomega_{i_1,a_1}\cdots\bomega_{i_k,a_k}\ \  {\rm{with}}\ \ a_p-a_{p-1}\in S_{i_p,i_{p-1},n},\  p\in[2,r].\end{equation}
The next result was proved in \cite{MY12a,MY12}. 
\begin{prop}\label{mysnake}  Suppose that $$\bomega=\bomega_{i_1,a_1}\cdots\bomega_{i_k,a_k},\ \ \bomega'=\bomega_{j_1,b_1}\cdots\bomega_{j_m,b_m}$$ satisfy \eqref{ps}, i.e. are prime snakes. 
\begin{enumerit} 
\item[(i)] If $b_1-a_k\in S_{j_1,i_k,n}$, then the $\hat\bu_n$--module $V(\bomega)\otimes V(\bomega')$ is reducible.
\item[(ii)] We have 
\begin{gather*}\wt_\ell V(\bomega)=\{\bomega(\underline p): \underline p\in   \mathbb P_{\bomegas}\},\ \ \dim V(\bomega)_\bpis=1 \ \ {\rm if}\ \  \bpi\in\wt_\ell V(\bomega),\\   \wt^+_\ell V(\bomega)= \{\bomega\}.\end{gather*} \item[(iii)] Suppose that $k=m$ and  $j_s=i_{s+1}$ and $b_s=a_{s+1}$ for $s\in[2,k-1] $. The following equality holds in $\cal K_0(\mathscr F_{n,\mathbb Z})$:
$$[V(\bomega)\otimes V(\bomega')] =[V(\bomega\bomega_{j_k,b_k})][V(\bomega'\bomega_{j_k}^{-1})]+[V(\bomega^+)][V(\bomega^-)],$$
 where $$\bomega^\pm = \prod_{p=1}^k\bomega_{\frac12(i_p+i_{p+1} \pm (a_p-a_{p+1})), \ \frac12(a_p+a_{p+1}\pm (i_p-i_{p+1}))},$$ and we understand $i_{k+1} = j_k$ and $a_{k+1}=b_k$. Moreover, 
$$\bomega^+\bomega^-\notin \wt_\ell(V(\bomega\bomega_{j_k,b_k})\otimes V(\bomega'\bomega_{j_k,b_k}^{-1})).$$

\item[(iv)] For $\underline p, \underline p'\in \mathbb P_\bomegas$ we have  
$$\bomega(p)\in\bomega(\underline p')
\cal Q^+  \iff p'(k)\geq p(k) \ \ {\rm{for\ all }}
 \ k\in [0,n+1].$$

\end{enumerit}
\hfill\qedsymbol\end{prop}

\subsubsection{} We note some consequences of Proposition \ref{mysnake}  for later use.
\begin{prop}\label{weights1} Let $i\in[1,n]$  and $\boa = (a_1,\cdots, a_r)$ with $a_1\le \cdots\le a_r$ for some $r\geq 1$ be such that $\bomega_{i,\boa}$ is a prime snake. For $j,s\in[1,r]$, with $j\le s$, set $\boa_{j,s}=(a_j,\cdots ,a_s)$ and $\boi_{j,s}=(i_j,\cdots, i_s)$.
\begin{enumerit}
\item[(i)] We have  $(p_1,\cdots, p_r)\in\mathbb P_{i,\boa}\implies (p_j,\cdots, p_s)\in\ \mathbb P_{i,\boa_{j,s}}.$
\item[(ii)] Conversely $(p_j,\cdots, p_s)\in\mathbb  P_{i,\boa_{j,s}}\implies (p_{i_1,a_1},\cdots, p_{i_{j-1},a_{j-1}}, p_j,\cdots, p_s, p_{i_{s+1},a_{s+1}}^*,\cdots , p_{i_r,a_r}^*)\in\mathbb P_{i,\boa}$.
 
\end{enumerit}
\end{prop}

\subsection{ A result of Kang, Kashiwara, Kim and Oh} We shall make crucial use of the main result of \cite{KKKO15} which we now recall.
 \begin{thm*}\cite[Theorem 3.12]{KKKO15} \label{simsoc} Suppose that $\bpi\in\cal P_n^+$ is such $V(\bpi)$ is real. Then for all $\bpi'\in\cal P_n^+$ the module $V(\bpi)\otimes V(\bpi')$ has simple head and simple socle.  Moreover the socle of $V(\bpi)\otimes V(\bpi')$ is the head of $V(\bpi')\otimes V(\bpi)$.\hfill\qedsymbol\end{thm*}

 \section{Higher order KR--modules and a prime factorization theorem}
 In this section we prove the analog of the main result of \cite{CP91} for elements  of  $\cal P_n^+$ concentrated at a node $i\in[1,n]$. We reformulate the notion  of a prime snake module $(i,\boa)$ in terms of  $(i,n)$--segments and define a notion of two $(i,n)$--segments being in general position. We show that   if $\bob\in\mathbb Z^r$ for some $r\ge 1$ then $V(\bomega_{i,\bob})$ can be written uniquely  (up to an overall permutation)  as a tensor product of modules associated to $(i,n)$--segments in general position.

\subsection{ The set $S_{i,n}$ and  $(i,n)$--segments}\label{segdef} 
For $i\in[1,n]$, let 
\begin{equation}\label{sin}
S_{i,n}=\{2j:1\le j\le \min\{i,n+1-i\}\}= S_{n+1-i,n}.\end{equation}
Notice that this is precisely the set $S_{i,i,n}$ defined in Section \ref{weyl}.
\subsubsection{Segments} \begin{defn} Say that an element $\boa=(a_1,\cdots, a_r)\in\mathbb Z^r$ is  an $(i,n)$-segment of length $r$ if   $a_p-a_{p-1}\in S_{i,n}$ for all $2\le p\le r$. Equivalently we say that $\boa$ is an $(i,n)$--segment if and only if $(i,\boa)$ is a prime snake of type $n$. \\\\
Set $$\boa^*= (a_1+n+1,\cdots, a_k+n+1),\ \ ^*\boa= (a_1-n-1,\cdots, a_k-n-1),$$ and notice that they are $(n+1-i,n)$--segments.
\hfill\qedsymbol \end{defn}
\noindent Since $0\notin S_{i,n}$ the entries of $\boa$ are all distinct and so in what follows we will  also think of segments as sets. \\\\
 {\bf Example.} If $n=3$ we have  $$S_{1,3}=\{2\}= S_{3,3},\ \ S_{2,3}=\{2,4\}.$$ The element $\boa=(0,4,6, 10)$ is 
 the union of the  three $(1,3)$-segments (and also $(3,3)$--segments) namely:  $(0)$, $(4,6)$ and $(10)$.  However $\boa$ is   a  $(2,3)$--segment of length 4.

\subsubsection{General and special position of segments}\label{genspec}  \begin{defn} Say that two $(i,n)$-segments $\boa=(a_1,\cdots, a_r)$ and $\bob=(b_1, \cdots, b_s)$  are in general position if their union does not contain an $(i,n)$--segment of length greater than $\max\{r,s\}$.  Otherwise we say that they are in special position. \hfill\qedsymbol\end{defn}
\noindent {\bf Examples.}
\begin{itemize}
\item[(i)] An $(i,n)$-segment is in general position with itself.
\item[(ii)] Consider the $(2,3)$--segments $\boa=(0,2,6,10), \ \ \bob=(4),\ \ \boc=(16, 18).$
Then $\boa$ and $\bob$ are  in special position since their union contains the $(2,3)$-segment $(0,2,4,6,10)$ while  \ $\boa ,\boc$ (and also $\bob$, $\boc$) are   in general position.\end{itemize} 
\subsection{ The KR--modules of type $(i,n)$ } 

Recall that for $i\in[1,n]$ and $\boa\in\mathbb Z^r$ we set$$\bomega_{i,\boa}=\bomega_{i,a_1}\cdots\bomega_{i,a_r}\in\cal P^+_{n}.$$ 
\begin{defn} Given  $\bomega\in\cal P^+_{n}$ we say that $V(\bomega)$ is a KR--module of type $(i,n)$ if there exists an $(i,n)$--segment $\boa$  such that $\bomega=\bomega_{i,\boa}$.\hfill\qedsymbol 
\end{defn}
\noindent We note that the  usual KR--module for $\hat\bu_n$ is of the form $\bomega_{i,\boa}$ where $\boa=(a,a+2,\cdots a+2r-2)$ for some $a\in\mathbb Z$ and $r\ge 1$. We  refer to the KR--modules associated with more general segments as the higher order KR--modules since they encode the reducibility data of  $V(\bomega_{i,a})\otimes V(\bomega_{i,b})$  in higher rank. 
\begin{rem}
The KR--modules of type $(i,n)$ are a special family of  the prime snake modules  \cite{MY12a,MY12}.
\end{rem}
\subsection{ A prime factorization result} We  state our first main theorem, which generalizes the result of  \cite{CP91}  in the rank one case. 
 \begin{thm}\label{main} 
 Let $i\in[1,n]$ and $\boa\in\mathbb Z^r$. There exists a unique integer $k\ge 1$ and  unique (up to a  permutation) $(i,n)$--segments $\boa_1,\cdots, \boa_k$   which are  in pairwise  general position such that $$V(\bomega_{i,\boa})\cong V(\bomega_{i,\boa_1})\otimes \cdots\otimes V(\bomega_{i,\boa_k}).$$ In particular $V(\bomega_{i,\boa})\otimes V(\bomega_{i,\boa})$  is irreducible  for all $\boa\in\mathbb Z^r$. In particular, $V(\bomega_{i,\boa})$ is prime if and only if $\boa$ is   an $(i,n)$--segment. \end{thm}
 The proof of the theorem occupies the rest of the section.
\subsection{Combinatorics of $(i,n)$--segments} We  give a more explicit formulation for a pair of  segments to be in special or general position. We use this  to prove that an element $\bomega_{i,\boa}$, $\boa\in\mathbb Z^r$  can be written uniquely (up to a renumbering) as a product  of elements associated to $(i,n)$--segments in general position.

\subsubsection{}   \begin{prop}\label{genpos} Let $r\ge m$ and assume that $\boa=(a_1,\cdots, a_r)$ and $\bob=(b_1,\cdots,\ b_m)$ are $(i,n)$--segments. 
\begin{enumerit}
\item[(i)]
The segments $\boa$ and $\bob$ are  in general  position if and only if one of the following holds:
\begin{enumerit}
\item[(a)]  $b_1-a_r>2\min\{i,n+1-i\}$ or   $a_1-b_m>2\min\{i,n+1-i\}$  
or  $b_1-a_1\notin 2\mathbb Z$,
\item[(b)]  $
\{b_1,\cdots, b_m\}\subset \{a_1,\cdots, a_r\} $.
\end{enumerit}
\item[(ii)]  The segments $\boa$ and $\bob$  are in special position if and only if there exists $1\le j\le m$  such that  one of the following hold:
\begin{enumerit}
\item[(a)]  $b_j-a_r\in S_{i,n}$
or $a_1-b_j
\in S_{i,n}$,
\item[(b)]  $b_1-a_1\in 2\mathbb Z$ and there exists $1\le k< r$ such that $a_k<b_j<a_{k+1}$.
\end{enumerit}
\end{enumerit}
\end{prop}
\begin{pf}  It  is clear that if  one of conditions in (i)(a) or if (i)(b) holds then $\boa$ and $\bob$ are    are in general position. For the converse  we suppose that  none of the conditions in $(a)$ and $(b)$ are satisfied and show that $\boa\cup\bob$ contains a segment of length $r+1$.  
If  $b_1>a_r$ (resp. $b_m<a_1$)  then $b_1-a_r\in S_{i,n}$ (resp. $a_1-b_m\in S_{i,n}$) which means that $(a_1,\cdots, a_r,b_1)$ (resp. $(b_m,a_1,\cdots, a_r)$) is an $(i,n)$-segment of length $r+1$ as needed. \\\\ Hence to complete the proof of  (i) we must consider the case when all of the following hold: $b_1\le a_r$, $ b_m\ge a_1$, and $ b_1-a_1\in 2\mathbb Z.$   If $ a_r< b_k$  (resp. $b_k<a_1$) for some $2\le k\le m$ (resp. $1\le k\le m-1$)  and $k$ is minimal (resp. maximal) with this property, then $b_{k-1}\le a_r<b_k$ (resp. $ b_k<a_1\le b_{k+1}$) and so  $(a_1,\cdots, a_r, b_k)$ (resp. $(b_k,a_1,\cdots a_r)$) is an  $(i,n)$-segment. Otherwise, we  have $a_1\le b_1 <b_2<\cdots <b_m\le a_r$.
Since condition $(b)$ does not hold, 
it follows that  $b_p\notin\{a_1,\cdots, a_r\}$ for some $1\le p\le m$; in other words,  there exists $2\le k\le r$ such that $a_{k-1}<b_p<a_k$.
 It follows that $(a_1,\cdots, a_{k-1}, b_k,a_k,\cdots a_r)$ is an $(i,n)$-segment. The proof of part (i) is complete. \\\\
 If either of the conditions in part (ii) hold, then  $\{a_1,\cdots a_r,b_j\}$ is an $(i,n)$--segment after applying a suitable permutation and hence $\boa$ and $\bob$ are in special position. Suppose that  $\boa$ and $\bob$ are in special position and assume that there does not exist $1\le j\le m$ and $1\le k< r$ with $a_k<b_j<a_{k+1}$. By part (i) we see that we cannot have $\{b_1,\cdots, b_m\}\subset \{a_1,\cdots, a_r\}$. Hence either there exists $j$ maximal with  $b_j<a_1$ or $j$ minimal with $b_j>a_r$. 
 In the first (resp. second) case either $j=m$ or $j<m$ and  $b_j<a_1\le b_{j+1}$ (resp. $j=1$ or $j>1$ and $b_{j-1}\le a_r<b_j$). If $j=m$ (resp. $j=1$) then 
 the segments can be in special position only if $a_1-b_m\in S_{i,n}$ (resp. $b_1-a_r\in S_{i,n}$). If $j<m$ (resp. $j>1$) then  $b_{j+1}-b_j\in S_{i,n}$ (resp. $b_j-b_{j-1}\in S_{i,n}$) and so $a_1-b_j\in S_{i,n}$ (resp. $b_j-a_r\in S_{i,n}$) and the proof is complete.
 \end{pf}
 
\begin{cor}
Suppose that $\boa=(a_1,\cdots, a_r)$ and $\bob=(b_1,\cdots, b_m)$ are $(i,n)$--segments in general position with $r\ge m$.  Then,
$$\{a_1,\cdots, a_r\}\cap\{b_1\cdots,b_m\}\ne\emptyset\implies \{b_1\cdots,b_m\}\subset\{a_1,\cdots, a_r\}.$$
\hfill\qedsymbol
\end{cor}

\subsubsection{} Given $\boa=(a_1,\cdots, a_r)\in \mathbb Z^r$ and $\bob=(b_1,\cdots, b_s)\in\mathbb Z^s$ set $$\boa\vee\bob=(a_1,\cdots, a_r, b_1,\cdots, b_s).$$
\begin{prop} \label{factorgenpos} Suppose that $\boa=(a_1,\cdots, a_r)\in\mathbb Z^r$. Then  $\bomega_{i,\boa}$  can be written uniquely (up to  a permutation) as a product $\bomega_{i,\boa}=\bomega_{i,\boa_1}\cdots\bomega_{i,\boa_k}$ where $\boa_1,\cdots,\boa_k$ are $(i,n)$--segments in pairwise general position.
\end{prop}
\noindent Before proving the proposition we give an example.\\\\
\noindent {\bf Example.}  Let  $\boa=(0,6, 4,2,10, 16, 10)\in\mathbb Z^7$, then the associated  $(2,3)$-segments  are $$\ \ \ \ \boa_1=(0,2,4,6,10), \ \ \ \boa_2=(10),\ \ \boa_3=(16),$$
or any permutation of these by an element of $S_3$

\begin{pf} Clearly it suffices to prove that there exists  a unique  integer $k\ge 1$ and  $(i,n)$--segments  $\boa_1,\cdots, \boa_k$   (which are unique up to a renumbering) in pairwise in general position and a permutation $\sigma\in S_r$ such that $$(a_{\sigma(1)},\cdots, a_{\sigma(r)})=\boa_1\vee\boa_2\vee\cdots\vee\boa_k.$$  We proceed by induction on $r$  with induction beginning trivially  at $r=1$. After applying an element of $S_r$ if needed, we may assume without loss of generality that:
$a_1\le a_s$ for all $1\le s\le r$ and that  $r_1\in[1,r]$  is maximal so that $\boa_1=(a_1, a_2,\cdots, a_{r_1})$ is an $(i,n)$-segment. If $r_1=r$ we are done  and otherwise we let $\bob=(a_{r_1+1},\cdots, a_r)$. 
The 
inductive hypothesis applies to $\bob$ and we let   $\boa_2,\cdots,\boa_k$ be the $(i,n)$--segments associated to $\bob$.\\\\
We  prove that $\boa_1$ is  in general position with $\boa_s$ for  $s\in[2,k]$. 
Suppose that $\boa_s=(b_1,\cdots, b_{m})$ and recall  that $a_1\le b_1$. Assume for a contradiction that $\boa_s$ and $\boa_1$ are in special position. If $m\le r_1$, then by Proposition \ref{genpos}(ii)
there exists $1\le p\le m$ such that either  $a_{k-1}<b_p<a_k$ for some $2\le k\le r_1$
or $a_{r_1}<b_p$ with $b_p-a_{r_1}\in S_{i,n}$. In either case after applying a suitable permutation if needed, we see that $(a_1,\cdots, a_{r_1}, b_p)$  defines
an $(i,n)$--segment  contradicting our choice of $r_1$. If $m>r_1$ then the maximality of $r_1$ implies that $b_1>a_1$ and also that $b_1-a_{r_1}\notin S_{i,n}$. Since by assumption $\boa_1$ and $\boa_s$ are in special position, this forces $b_1\le a_{r_1}$. Choose $p\in[1,r_1]$ minimal so that  $a_p<b_1\le a_{p+1}$. Then,   $(a_1,\cdots, a_p,b_1,\cdots, b_m)$ is a $(i,n)$--segment
 and again we have a contradiction to the maximality of $r_1$. This proves that $\boa_1$ and $\boa_s$ are in general position. 
\\\\
It remains to prove that $k$ is unique and that the  segments are unique up to an element of $S_k$. For this, suppose that $\boc_1,\cdots,\boc_\ell$ is another set of $(i,n)$--segments in pairwise  general position and assume that  $(a_{\sigma(1)},\cdots, a_{\sigma(r)})=\boc_1\vee\cdots\vee \boc_\ell$ for some $\sigma\in S_r$. Since $a_1$ is minimal it must occur as the first term in $\boc_p$ for some $1\le p\le \ell$ and   assume without loss of generality that $p=1$ and also that $\boc_1$   has maximal length say $s_1$ amongst those $\boc_s$ with first term $a_1$.  Since $r_1$ is the maximum length of an $(i,n)$-segment starting at $a_1$ we have $s_1\le r_1$. We claim that $\boc_1=\boa_1$. Otherwise,  there exists $1<p\le r_1$ minimal such that $a_p$ does not occur in $\boc_1$.    All the other segments $\boc_s$ whose initial term is  $a_1$ have length at most $s_1$ and hence by Corollary \ref{genpos} must be contained in $\boc_1$.  Hence none of these segments contain $a_p$  and so  there must exist an $(i,n)$--segment $\boc_j$ of length $s_j$ containing $a_p$. The minimal term of  $\boc_j$ is $a_m$ for some $m$ satisfying  $a_p\ge a_m>a_1$. This means that $a_{m-1},a_m\in\boa_1$ and so $a_m-a_{m-1}\in S_{i,n}$. Consider $\boc_1\cup\boc_j$. If $s_1\ge s_j$ then since $a_{p-1}\in\boc_1$ it follows that $\boc_1\cup\{a_p\}$ is a longer $(i,n)$-segment in the union.  If  $s_j>s_1$   the preceding discussion shows that $\{a_{m-1}\}\cup \boc_j$ is a longer segment in the union. In any case we have a contradiction to the fact that $\boc_1$ and $\boc_j$ are in general position. It follows that $\boa_1=\boc_1$ and that $\bob=\boc_2\vee\cdots\vee\boc_\ell$. The uniqueness is now immediate by the inductive hypothesis.

\end{pf}

 \subsubsection{} We need a further equivalent formulation on segments in special position. We remark that in the language of \cite{MY12} the next proposition proves that any pair of $(i,n)$--segments in special position contains a pair of segments which form an extended $T$--system. In the following,  given any $(i,n)$--segment $\boc=(c_1,\cdots, c_k)$ we adopt the convention that  $c_m=-\infty$ if $m<1$ and $c_m=\infty$ if $m>k$.  
 \begin{lem}\label{tsys} Let  $\boa=(a_1,\cdots, a_r)$ and $\bob=(b_1,\cdots, b_s)$ be $(i,n)$--segments  with $r\ge s$. Then $\boa$ and $\bob$ are in special position if and only if $b_1-a_1\in 2\mathbb Z$ and there exists $p\in \mathbb Z_{\geq 0}$,  $\boa_1=(a_j,\cdots, a_{j+p})\subset \boa$ and $\bob_1=(b_m,\cdots, b_{m+p})\subset\bob$ such that \begin{equation}\label{minmax} \max\{ a_{j-1}, b_{m-1}\}<\min\{a_j,b_m\},\ \ \min\{a_{j+p+1}, b_{m+p+1}\}>\max\{a_{j+p}, b_{m+p}\},\end{equation}  and  either\begin{equation}\label{ets} a_{j+k}=b_{m+k-1}\ \ \ {\rm{for \ all}}\ k\in[1,p] \ \ {\rm{or}}\ \ a_{j+k-1}=b_{m+k}\  {\rm{for \ all}}\ k\in[1,p].\end{equation}  
 \end{lem}
 Before proving the Lemma we give an example.\\\\
 {\bf Example.}
   Consider the  $(2,3)$--segments $$\boa = (-4,0,4,8,10,14),\ \ \bob = (0,2,4,8,12,14).$$ They  are in special position since their union contains the $(2,3)$--segment
    $ (-4,0,2,4,8,10,12,14)$. \\ 
    In the notation of the lemma, taking $p=2$, $j=3$ and $m=2$ we have $\boa_1 = (4,8,10)$, 
    $\bob_1 = (2,4,8)$, which clearly satisfy  
    $$\max\{a_2,b_1\}=0<2=\min\{a_3,b_2\}, \ \ \min\{a_6,b_5\} = 12> 10=\max\{a_5,b_4\}$$ and \eqref{ets} as well.

 \begin{pf} Suppose that  $\boa$ and $\bob$ satisfy  \eqref{minmax}  and that $a_{j+k}=b_{m+k-1}$ for $k\in[1,p]$ (resp. $a_{j+k-1}=b_{m+p}$) for all $k\in[1,p]$. Then  $\{a_1,\cdots, a_{j+p}, b_{m+p}, a_{j+p+1},\cdots, a_r\}$ is an $(i,n)$--segment of length $r+1$ which proves the converse direction.\\\\
 For the forward direction, 
 we proceed by induction on $s$.
 To see that induction begins at  $s=1$ notice that if $r=1$ there is nothing to prove.   If $r>1$ then by Proposition \ref{genpos} we have  $a_j<b_1<a_{j+1}$ for some $0\le j\le r$.   Clearly \eqref{minmax} and \eqref{ets} hold by taking  $p=0$,  $\boa_1=(a_{j+\delta_{j,0}})$ and $  \bob_1=(b_1)$.
  For  the inductive step  assume that the result holds   for all $(i,n)$--segments in special position for all  $s'<s$ and  all $r\ge s'$. We  consider two main cases.\\\\
 For the first case we assume   that $b_1\leq a_1$.  If $b_2<a_1$ then  observe that $\boa$ and $\bob'=(b_2,\cdots b_s)$ are in special position. Hence the inductive hypothesis applies to these pairs and so we can choose $\boa_1$ and $\bob_1$ so that they satisfy \eqref{minmax} and \eqref{ets}. If
  $b_1=a_1$
 or  $b_1< a_1$ and $a_k<b_j<a_{k+1}$ for some $j\in[2,s]$ and $k\in[1,r]$ then the pairs 
 $\boa'=(a_2,\cdots, a_r)$ and $\bob'=(b_2,\cdots, b_s)$   are in special position.  The induction hypothesis applies to these pairs, and so we can choose  $\boa_1$ and $\bob_1$ satisfying \eqref{minmax} and \eqref{ets}. To complete the first case we  must consider the situation when
  $b_1<a_1\le b_2$ 
 and  there does not exist $j\ge 2$ with $a_k<b_j<a_{k+1}$.  Since $r\ge s$ there exists  $p\in[2,s+1]$ maximal so that $a_{m-1}=b_m$ if $2\le m\le p-1$ and $a_{p-1}\neq b_p$. In particular, this means that $a_{p-1}<a_p\le b_p$ and  so the elements $\bob_1=(b_1,\cdots, b_{p-1})$ and $\boa_1=(a_1,\cdots,a_{p-1})$ satisfy the conditions of the proposition and give the inductive step.\\\\
For the second case we assume that  $b_1>a_1$. Working in a similar fashion from the other side we see that we can further reduce to the case when  $b_s<a_r$ as well. 
Choose $p$ maximal and $p'$ minimal so that $a_p<b_1<b_s<a_{p'}$. Since $\boa$ and $\bob$ are in special position, Proposition \ref{genpos} shows that there exists $j\in [1,s]$ minimal such that $a_k<b_j<a_{k+1}$ for some $k\in [p,p']$. Define an integer $m$ as follows: if either $j=1$ or $j>1$ and $b_{j-1}<a_k$ we take $m=0$; otherwise we take $m$ so that $a_{k-m'}=b_{j-m'-1}$ for all $0\leq m'<m$ and $a_{k-m}> b_{j-m-1}$. Notice that $m$ must exist since $a_p<b_1$ and $j$ was chosen minimal. This time we take
$$ \boa_1 = (a_{k-m},\cdots, a_k), \ \   \bob_1 = (b_{j-m},\cdots, b_j).$$
 In particular we have 
 $$a_{k-m-1}, b_{k-m-1}< a_{k-m}<b_{k-m}, \ \ \ a_k<b_j<\min\{a_{k+1},b_{k+1}\},$$
  and the proof of the inductive step is complete. 

 \end{pf}

\subsection{Proof of Theorem \ref{main}: Reducibility}\label{redpf}

Assume that 
 $\boa=(a_1,\cdots, a_r)$ and $\bob=(b_1,\cdots, b_s)$ are $(i,n)$--segments in special position with  $r\ge s$.  We prove that $V(\bomega_{i,\boa})\otimes V(\bomega_{i,\bob})$ is reducible by showing that  \begin{equation}\label{strict}\wt_\ell \left(V(\bomega_{i,\boa})\otimes V(\bomega_{i,\bob})\right)\ne \wt_\ell V(\bomega_{i,\boa}\bomega_{i,\bob}).\end{equation} 
Choosing $\boa_1$ and $\bob_1$ as in Lemma \ref{tsys}, we use Proposition \ref{mysnake} to choose elements $\bpi$ and $\bomega$ in $\cal P_n$ as follows: \begin{equation}\label{omega}\bpi\in \wt_\ell \left(V(\bomega_{i,\boa_1})\otimes V(\bomega_{i,\bob_1})\right)\setminus \wt_\ell V(\bomega_{i,\boa_1}\bomega_{i,\bob_1}),\ \ \bomega=\bomega_{i,\boa_0}\bomega_{i,\bob_0}\bpi\bomega_{n+1-i, \boa_2^*}^{-1}\bomega_{n+1-i, \bob_2^*}^{-1}.\end{equation} Here $\boa_0=(a_1,\cdots, a_{j-1})$, $\boa_2=(a_{j+p+1},\cdots, a_r)$ and $\bob_0$, $\bob_2$ are defined similarly, and we  understand that these segments can be empty. Using \eqref{minmax} we see that  $\bpi$ is in the subgroup of $\cal P_n^+$  generated by elements $\bomega_{j,c}$ with $$\min\{a_{j-1},b_{m-1}\}<c\leq  \max\{a_{j+p}, b_{m+p}\}+n+1<\min\{a_{j+p+1}, b_{m+p+1}\}+n+1.$$ It follows that  the expression  in \eqref{omega} for $\bomega$ is reduced (c.f. Section \ref{snakepaths}).
Writing $\bpi=\bomega_1\bomega_2$ with $\bomega_1\in \wt_\ell V(\bomega_{i,\boa_1})$ and $\bomega_2\in V(\bomega_{i,\bob_1})$ we see by using Proposition \ref{weights1} that $$\bomega=(\bomega_{i,\boa_0}\bomega_1\bomega_{i+1-i,\boa_2^*}^{-1})(\bomega_{i,\bob_0}\bomega_2\bomega_{n+1-i,\bob_2^*}^{-1})\in \wt_\ell V(\bomega_{i,\boa})\otimes V(\bomega_{i,\bob}).$$  The assertion in  \eqref{strict} follows if we prove that \begin{equation}\label{notintrip} \bomega\notin\wt_{\ell} \left(W(\bomega_{i,\boa_0}\bomega_{i,\bob_0})\otimes V(\bomega_{i,\boa_1}\bomega_{i,\bob_1})\otimes W(\bomega_{i,\boa_2}\bomega_{i,\bob_2})\right).\end{equation}
 This is because  the module $V(\bomega_{i,\boa}\bomega_{i,\bob})$ occurs in the Jordan--Holder series of the triple tensor product. Suppose for a contradiction that $$\bpi_1\bpi_2\bpi_3=\bomega= \bomega_{i,\boa_0}\bomega_{i,\bob_0}\bpi\bomega_{n+1-i,\boa_2^*}^{-1}\bomega_{n+1-i,\bob_2^*}^{-1}$$ where  $\bpi_1$, $\bpi_2$, $\bpi_3$ are $\ell$--weights of the corresponding modules in the tensor product. This gives, $$\bpi_1\bomega_{i,\boa_0}^{-1}\bomega_{i,\bob_0}^{-1}=\bpi_2^{-1}\bpi_3^{-1}\bpi\bomega_{n+1-i,\boa_2^*}^{-1}\bomega_{n+1-i,\bob_2^*}^{-1}.$$ If $\bpi_1\ne \bomega_{i,\boa_0}\bomega_{i,\bob_0}$ then some $\bomega_{i,c}^{-1}$ with $c\in\boa_0\cup\bob_0$ survives in a reduced expression of the left hand side. On the other hand \eqref{minmax} shows that the right hand side is 
  in the subgroup of $\cal P_n^+$ generated by $\bomega_{j,d}$ with $d>\max\{a_{j-1}, b_{m-1}\}$ and hence cannot have $\bomega_{i,c}^{-1}$ in a reduced expression and  we have a contradiction. It follows that $\bpi_1=\bomega_{i,\boa_0}\bomega_{i,\bob_0}$.
   Similarly as before write
  $$\bpi_3\bomega_{n+1-i,\boa_2* }\bomega_{n+1-i,\bob_2^*} = \bpi \bpi_2^{-1}.$$
  If $\bpi_3\neq \bomega_{n+1-i,\boa_2^*}^{-1}\bomega_{n+1-i,\bob_2^*}^{-1}$ there exists $c\in\boa_2^*\cup \bob_2^*$ such that $\bomega_{i,c}$ survives in a reduced expression of the left hand side. Another application of  \eqref{minmax} shows that $\bpi_2\bpi$ is in the subgroup of $\cal P_n^+$ generated by $\bomega_{j,d}$, $$d\leq \max\{a_{j+p},b_{m+p}\}+n+1<\min \{a_{j+p+1},b_{m+p+1}\}+n+1\leq c,$$
  and hence we a have a contradiction. Then we must also have $\bpi_3 = \bomega_{n+1-i,\boa_2^*}^{-1}\bomega_{n+1-i,\bob_2^*}^{-1}$ and thus $\bpi = \bpi_2$. But  this again contradicts our choice of $\bpi$ and so the assertion in \eqref{notintrip} holds and the proof of reducibility is complete.

\subsection{Proof of Theorem \ref{main}: Irreducibility}\label{irredpf}

We prove  irreducibility  by induction on $k$.  Let $k=2$ and $\boa$ and $\bob$ be $(i,n)$--segments in general position. The irreducibility of $V(\bomega_{i,\boa})\otimes V(\bomega_{i,\bob})$ follows from the fact that $\dim (V(\bomega_{i,\boa})\otimes V(\bomega_{i,\bob}))_{\bomegas_{i,\boa}\bomegas_{i,\bob}}=1$ and the following claim: \begin{equation}\label{sochda}\soc(V(\bomega_{i,\boa})\otimes V(\bomega_{i,\bob}))\cong V(\bomega_{i,\boa}\bomega_{i,\bob})\cong{\rm{hd}}(V(\bomega_{i,\boa})\otimes V(\bomega_{i,\bob})).\end{equation} For  the first isomorphism assume that  $V(\bpi)$ is in the socle of $V(\bomega_{i,\boa})\otimes V(\bomega_{i,\bob})$ for some $\bpi\in\cal P_n^+$. 
Using Lemma \ref{socle}, write \begin{equation}\label{pisoc}\bpi=\bomega_{i,\boa}\bomega(\underline p),\ \ 
(\bpi^*)^{-1}\bomega_{n+1-i,\bob^*}=\bomega(\underline g),\ \ \underline p\in\mathbb P_{i,\bob},\ \ \underline g\in\mathbb P_{i,\boa}.\end{equation} Since $\bpi\in\cal P_n^+$ it follows from \eqref{pjm} of Section \ref{gjm} that there exist integers $0\le m_1,\cdots,  m_s\le \min\{i, n+1-i\}$ so that $\underline p=(p_{i,b_1}^{m_1},\cdots ,p_{i,b_s}^{m_s})$.\\\\
We use the equivalent formulation  given in Proposition \ref{genpos} for a pair of $(i,n)$--segments to be  in general position. Without loss of generality we can assume that either $a_1-b_s>2\min\{i,n+1-i\}$ or that $a_1-b_1\notin 2\mathbb Z$ or $\{b_1,\cdots, b_s\}\subset\{a_1,\cdots, a_r\}$. \\\\
Suppose that  $a_1-b_s>2\min\{i,n+1-i\}$ or that $a_1-b_1\notin 2\mathbb Z$.  Then $b_j+2m_j\notin\{a_1,\cdots, a_r\} $ and since   $\bpi\in\cal P_n^+$ we must have  $m_j=0$ for all $j\in[1,s]$ proving $\bpi=\bomega_{i,\boa}\bomega_{i,\bob}$.\\\\
Consider the case when $\{b_1,\cdots , b_s\}\subset\{a_1,\cdots, a_r\}$. Suppose for a contradiction that there exists $k$ minimal with $m_k\ne 0$ and choose $j\in[1,r]$ with $b_k=a_j$. The first equality in \eqref{pisoc} can be written as:$$\bpi=\bomega_{i, a_1}^{\epsilon_1}\cdots\bomega_{i,a_{j-1}}^{\epsilon_{j-1}}\bomega_{i,a_j}\bpi_1,\ \ \bpi_1=\bomega_{i,a_{j+1}}\cdots\bomega_{i,a_r}\bomega(p_k)\cdots\bomega(p_r)\in\cal P_n^+,$$ where  $\epsilon_m=1$ if $a_m\notin\{b_1,\cdots,b_{k-1}\}$ and $\epsilon_m=2$ otherwise. 
Substituting in the second equality in \eqref{pisoc}, we get $$\bomega_{n+1-i, a_1+n+1}^{-1}\cdots\bomega_{n+1-i,a_{j-1}+n+1}^{-1}(\bpi_1^*)^{-1}\bomega_{n+1-i,b_{k+1}+n+1}\cdots\bomega_{n+1-i,b_s}=\bomega(\underline g). $$ 
Notice that  $\bpi_1$ is in the submonoid  of $\cal P_n^+$ generated by elements $\bomega_{p,c}$ with $c> a_j$. Hence $\bomega_{n+1-i,a_j+n+1}^{-1}$ does not occur on the left hand side of the preceding equation and so cannot occur on the right. This means that $\bomega(g_j)$ has an expression of the form $\bomega_{p,c}$ with $c<a_j+n+1$  which does not occur on the left hand side of the previous equality. This gives the desired contradiction and so $\bomega(\underline p)=\bomega_{i,\bob}$ and the first isomorphism of \eqref{sochda} is proved.\\\\
If $\boa=\bob$ then the second isomorphism follows by duality. In particular this proves that the module $V(\bomega_{i,\boa})$ is real.
If $\boa\ne\bob$ then using Theorem \ref{simsoc}  the second isomorphism follows if we prove that 
$$\soc (V(\bomega_{i,\bob})\otimes V(\bomega_{i,\boa}))\cong V(\bomega_{i,\boa}\bomega_{i,\bob}).$$
Assume that $V(\bpi)$ is in the socle and write $$\bpi=\bomega_{i,\bob}\bomega(\underline p),\ \ 
\underline p=(p_1,\cdots, p_r)\in\mathbb P_{i,\boa}.
$$
By \eqref{pjm} there exist integers $m_j\in[0,\min\{i,n+1-i\}]$ such that $p_j=p_{i,a_j}^{m_j}$. If  $a_1-b_s>2\min\{i,n+1-i\}$ or if $b_1-a_1\notin 2\mathbb Z$ then $a_j+2m_j\notin \{b_1,\cdots, b_s\}$ for $m_j> 0$.  Hence  the restriction that $\bpi\in\cal P_n^+$ forces $m_j=0$ for all $j\in[1,r]$ and $\bpi=\bomega_{i,\bob}\bomega_{i,\boa}$. 
\\\\
Suppose that  $\{b_1,\cdots, b_s\}=\{a_{i_1},\cdots, a_{i_s}\}$ and  that there exists $j\in[1,r]$ maximal with the property that  $m_j>0$. Since $g_j(i)=a_j+2m_j$ the condition that $\bpi\in\cal P_n^+$ forces $a_j+2m_j=a_{i_p}$  for some $p\in[1,s]$ with $i_p>j$. But the maximality of $j$ implies that $g_{i_p}(i)=a_{i_p}$ and we get a contradiction to $\underline g\in\mathbb P_{i,\boa}$. Hence we have proved the second isomorphism in \eqref{pisoc} and the irreducibility is established when $k=2$.\\

For the inductive step, assume that $\boa_1,\cdots,\boa_k$ are $(i,n)$--segments in general position of length $r_1,\cdots, r_k$ and that the irreducibility  is known for all $k'<k$. To prove it for $k$, we proceed by a further induction on $N=r_1+\cdots+r_k$ . This induction clearly begins at $N=k$ by Proposition \ref{weyl} (iii). Let $a$ be the maximal element of the set $\boa_1\cup\cdots\cup\boa_k$ and let $r$ be the number of segments with $a=\max \boa_s$ and assume without loss of generality that $\boa_1,\cdots,\boa_r$ are the segments whose maximal value is $a$. Let $\boa_s'=\boa_s\setminus\{a\}$, $1\leq s\leq k$. Then $\boa_1',\cdots, \boa_k'$ are $(i,n)$--segments in general position and the sum of their lengths  $N-r$.  The inductive hypothesis implies that $V(\bomega_{i,\boa_1'})\otimes\cdots \otimes V(\bomega_{i,\boa_k'})$ is irreducible and hence a quotient of $W(\bomega_{i,\boa_1'}\cdots\bomega_{i,\boa_k'})$.  By Proposition \ref{weyl}(iii) we have $$V(\bomega_{i,a})^{\otimes r}\otimes W(\bomega_{i,\boa_1'}\cdots\bomega_{i,\boa_k'})\cong W(\bomega_{i,\boa_1}\cdots\bomega_{i,\boa_k}).$$ It follows that
the quotient  $V(\bomega_{i,a})^{\otimes r}\otimes V(\bomega_{i,\boa_1'})\otimes \cdots\otimes V(\bomega_{i,\boa_k'})$ is $\ell$--highest weight and hence so is its quotient $V(\bomega_{i,a})^{\otimes  (r-1)}\otimes V(\bomega_{i,\boa_1})\otimes  V(\bomega_{i,\boa_2'})\otimes \cdots\otimes V(\bomega_{i,\boa_k'})$. Again using the fact that we have shown that  $V(\bomega_{i,a})\otimes V(\bomega_{i,\boa_1})$ is irreducible it follows that $V(\bomega_{i,\boa_1})\otimes  V(\bomega_{i,a})^{\otimes  (r-1)}\otimes  V(\bomega_{i,\boa_2'})\otimes \cdots\otimes V(\bomega_{i,\boa_k'})$ is $\ell$--highest weight. An iteration of this argument proves that $V(\bomega_{i,\boa_1})\otimes\cdots\otimes V(\bomega_{i,\boa_k})$ is $\ell$--highest weight. Working with the duals as usual establishes irreducibility.

\section{Inflation and Monoidal Categorification}\label{secinflation} 
Throughout this section we work with a triple of positive integers $(i,\bar i, n)$ which are related by requiring \begin{equation}\label{triple} n+1=i(\bar i+1).\end{equation}
Equivalently if we fix $n$ then $i$ can vary over all divisors of $n+1$ and $\bar i$ is defined by the preceding equation.
We  study the relationship between higher order Kirillov--Reshetikhin modules for $\hat{\bu}_{\bar i}$ and $\hat\bu_n$. We also study  connections of our theory 
 with the work of \cite{HL10, HL13} on monoidal categorification of cluster algebras.
\subsection{The homomorphism $\Phi_{\bar i,n}$ and the category $\mathscr C_{\bar i,n}$}\label{Phi}
Let $\Phi_{\bar i,n}:\cal P_{\bar i}\to\cal P_n$ be the group homomorphism defined by extending, $$\Phi_{\bar i,n}(\bomega_{j,a})=\bomega_{ij,ia},\ \ j\in[1,\bar i], \ \ a\in\mathbb Z.$$ We shall continue to denote by $\Phi_{\bar i, n}$ the induced homomorphism $\mathbb Z[\cal P_{\bar i}]\to\mathbb Z[\cal P_n]$. Notice that for $j_1,j_2\in[1,\bar i]$ and $a,b\in\mathbb Z$ we have,
$$ a-b\in S_{j_1,j_1,\bar i}\implies i(a-b)\in S_{ij_1,ij_2,n}.$$ It follows that if 
$\bomega\in\cal P_{\bar i}^+$ is a prime snake then 
 $\Phi_{\bar i, n}(\bomega)$ is a prime snake in $\cal P_n^+$. \\\\ 
Let $\cal H_{\bar i, n}$ (resp. $\cal H_{\bar i,n}^+$) be the subgroup (resp. submonoid) of $\cal P_{n}$ generated by the elements of the set $$\{\bomega_{ij, ia}:  j\in[1,\bar i], \ 
j-a\in 2\mathbb Z,\ \  a\in(-\infty ,0]\}.$$
 Let  $\mathscr C_{\bar i,n}$ be the full subcategory of $\mathscr F_{n}$ consisting of objects whose Jordan--Holder components are of the form $V(\bomega)$, $\bomega\in\cal H_{\bar i,n}^+$. 
 Notice that  $\mathscr C_{n,n} $ is precisely the category $\mathscr C^-$ studied in \cite{HL13}.  Let $\mathscr C_{\bar i, n}^1$ be the full subcategory consisting of objects in $\mathscr C_{\bar i, n}$  whose objects have Jordan--Holder constituents $V(\bpi)$ where $\bpi$ in the subgroup  $\cal H_{\bar i, n}^1$ generated by $\bomega_{ij,ic}\in\cal H_{\bar i, n}$, with $c\ge -3$. This is the analog of the category $\mathscr C_1$  studied in \cite{HL10}.\\
\subsection{}   Recall from  equation \eqref{truncg} the definition of the truncation of the $q$--character of an object of $\mathscr F_n$ at a subgroup $\cal G$ of $\cal P_n$.  Given $\boa = (a_1,\cdots, a_r)\in \mathbb Z^r$ it will be convenient to set $$i\boa = (ia_1,\cdots, ia_r).$$
 The following is  the main result  of this section.
 \begin{thm}\label{infliso}
 \begin{enumerit}
\item[(i)] The category $\mathscr C_{\bar i,n}$ is a monoidal category and the Grothendieck ring $\cal K_0(\mathscr C_{\bar i,n})$ is a polynomial algebra generated by the elements $\{ [V(\bomega_{ij,ia})]:\bomega_{ij,ia}\in\cal H_{\bar i,n}\}$. 
\item[(ii)] The assignment $$\chi^{\cal H_{\bar i,n}}: K_0(\mathscr C_{\bar i,n})\to \mathbb Z[\cal H_{\bar i,n}],$$ is an injective homomorphism and  \begin{equation}\label{phichi}\Phi_{\bar i,n}\circ\chi^{\cal H_{\bar i,\bar i}}([V(\bomega_{j,a})])= \chi^{\cal H_{\bar i,n}}([V(\Phi_{\bar i,n}(\bomega_{j,a}))]),\ \ \bomega_{j,a}\in\cal H_{\bar i,\bar i}^+.\end{equation}
\item[(iii)] We have an isomorphism $\Psi_{\bar i, n}:\cal K_0(\mathscr C_{\bar i,\bar i})\to \cal K_0(\mathscr C_{\bar i,n})$ such that  $$\Psi_{\bar i, n}([V(\bomega)])=[V(\Phi_{\bar i,n}(\bomega))],$$ where $\bomega\in\cal H_{\bar i,\bar i}$ is  of snake type or $\bomega=\bomega_{i,\boa}$ for some $\boa\in \mathbb Z^r$.\\
\item[(iv)] The restriction of $\Psi_{\bar i,n}$ to $\cal K_0(\mathscr C_{\bar i,\bar i}^1)$ satisfies 
\begin{equation}\label{primeinf}\Psi_{\bar i,n}([V(\bomega)])=[V(\Phi_{\bar i,n}(\bomega))],\ \ \bomega\in\cal H_{\bar i,\bar i}^1.\end{equation} In particular $\mathscr C_{\bar i, n}^1$ is a monoidal tensor category and  $\cal K_0(\mathscr C_{\bar i,n}^1)$ is a monoidal categorification of a cluster algebra of type $A_{\bar i}$. 
\end{enumerit}
 \end{thm}
\begin{rem} Parts (i) and (ii) of the theorem were proved for $\mathscr C_{\bar i, \bar i}$ in \cite{HL13}.  For a more detailed description of the relationship with cluster algebras we refer the reader to \cite[Section 4.3]{HL10} and \cite[]{BC19a}.
\end{rem}
\noindent{\bf Conjecture. }  We conjecture that $$\Psi_{\bar i,n}([V(\bomega)])=[V(\Phi_{\bar i,n}(\bomega))],\ \ \bomega\in\cal H_{\bar i,\bar i}.$$ Clearly parts (iii) and (iv) of the theorem establish special cases of the conjecture.\\\\
The proof of the theorem is in several steps.\\\\
{\em From now on for ease of notation we shall set $\Phi_{\bar i, n}=\Phi$ and $\Psi_{\bar i, n}=\Psi$.} 

\subsection{Properties of $\Phi$}
\subsubsection{} \label{inflatedlattice}   Let  $\phi: P_{\bar i}\to P_n$ be the group  homomorphism defined  by extending the assignment $  \phi(\omega_j)= \omega_{ij},\ \ j\in[1,\bar i]$.   Clearly $\wt\circ\ \Phi=\phi\circ\wt. $
\begin{lem}\label{partialorder}
We have $\phi(Q_{\bar i}^+)\subset Q_n^+$. Moreover for $\bpi,\bomega\in\cal P_{\bar i}$, $$\wt\bomega-\wt\bpi\in Q_{\bar i}^+\iff \phi(\wt\bomega)-\phi(\wt\bpi)\in Q_n^+.$$
\end{lem}
\begin{pf} Writing $\alpha_j=2\omega_j-\omega_{j-1}-\omega_{j+1}$, it is easily checked that $$\phi(\alpha_j)=\sum_{p=1}^{i}
\sum_{s=p-i}^{i-p}\alpha_{ij+s},\ \ j\in[1,\bar i].$$
The forward direction of the second assertion of the lemma is  now immediate. For the converse assume that  $\wt\bomega-\wt\bpi=\sum_{j=1}^{\bar i} s_j\alpha_j$ with $s_m<0$ for some $m\in[1,\bar i]$. The assertion follows once we notice  that the coefficient of $\alpha_{im}$  in $(\phi(\wt\bomega)-\phi(\wt\bpi))$ is  $is_m$.

\end{pf}
\begin{rem}\label{infroot} More generally, one computes that   $\Phi_{\bar i,n}(\balpha_{j,a})$, $j\in [1,\bar i]$, $a\in \mathbb Z$ is equal to
$$ \left(\prod_{k=1}^{i-1} \prod_{p=1}^k \balpha_{i(j-1)+k,\, i(a+1) -k+2p -2}\balpha_{i(j+1)-k,\, i(a+1) -k+2p -2}\right) \prod_{p=1}^i \balpha_{ij,\, i(a+1) -i+2p -2}.$$\\
\end{rem}

\subsubsection{The set $\Phi(\wt_\ell V(\bomega_{k,a}))$}\label{wtphi}  Given $a\in\mathbb Z$ and $k\in[1,\bar i ]$  we would like to show that $\Phi(\bomega)\in \wt_\ell V(\bomega_{ki,ai})$ for each $\bomega\in \wt_\ell V(\bomega_{k,a})$.  Since $\bomega= \bomega(p)$ for some $p\in \mathbb P_{k,a}$, it becomes natural in view  of the results of Section \ref{MYsec}, to define an element $\Phi(p)\in \mathbb P_{ki,ai}$ such that $\bomega(\Phi(p))=\Phi(\bomega(p))$. \\\\
For  $p\in \mathbb P_{k,a}$  define $\Phi(p):[0,n+1]\to \mathbb Z$ as follows:   
 for $j\in[0,\bar i+1]$ and $ j'\in[0,i-1]$ with $ij+j'\in[0,n+1]$, $$ \Phi(p)(ij+j')=\begin{cases} ip(j) +j', \ \  p(j+1)-p(j) =1,\\ ip(j)-j',\ \   p(j+1)-p(j)=-1.
 \end{cases}$$  A straightforward checking shows that $$\Phi(p)\in\mathbb P_{ik,ia}\ \ {\rm{and}}\ \  \bomega(\Phi(p))=\Phi(\bomega(p)).$$ 
Conversely, suppose that $g\in\mathbb P_{ik,ia}$ is such that $\bomega(g)\in \Phi( \cal P_{\bar i})$; equivalently:$$m\in\boc_g^-\cup\boc_g^+\cup\{0,n+1\}\implies m,\  g(m)\in i\mathbb Z.$$ 
 Let $r\in[0,\bar i+1]$ and  choose $m\in\boc_g^-\cup\boc_g^+\cup\{0,n+1\}$ with $$m\le ir\ \ {\rm{and}}\ \   g(ir)=g(m)\pm (ir-m).$$ It is immediate that $g(ir)$ is also in $i\mathbb Z$ and that $g(ir+j)=g(ir)\pm j$ for all $1\le j\le i$.
Define $$\Phi^{-1}(g):[0,\bar i+1]\to\mathbb Z,\ \ \Phi^{-1}(g)(r)=g(ir)/i,\ \ r\in[0,\bar i+1].$$
  It is straightforward to see that $\Phi^{-1}(g)\in\mathbb P_{k,a}$ and that  $$\bomega(g)=\Phi(\bomega(\Phi^{-1}(g))). $$ 
  The following is trivial.
  \begin{lem} Suppose that $s\in\{1,2\}$ and $p_s\in\mathbb P_{k_s, a_s}$ for some $p_s\in[1,\bar i]$ and $a_s\in\mathbb Z$.
  Then $$p_1(j)<p_2(j)\ {\rm{for \ all }}\ j\in[0,\bar i+1]\implies \Phi(p_1)(j)<\Phi(p_2)(j)\  {\rm{for \ all }}\ j\in[0,n+1].$$ A similar assertion holds for $s=1,2$ and  $g_s\in\mathbb P_{ik_s,ia_s}$ satisfying $\bomega(g_s)\in\Phi(\cal P_{\bar i})$.\hfill\qedsymbol
  \end{lem}
\subsubsection{The set $\wt_\ell V(\Phi(\bomega))$} The following proposition is crucial in the proof of Theorem \ref{infliso}.
\begin{prop}\label{inflpaths}
   Suppose that $\bomega\in\cal P_{\bar i}^+$ is a prime snake. For $\bpi\in \cal P_{\bar i}$ we have
$$\bpi\in \wt_\ell V(\bomega)\iff \Phi(\bpi)\in\wt_\ell V(\Phi(\bomega)) .$$
\end{prop}
\begin{pf}
Writing $\bpi=\bomega(p_1)\cdots\bomega(p_r)$ with $(p_1,\cdots, p_r)\in\mathbb P_{\bomegas}$, it is immediate from Lemma  \ref{wtphi} that $(\Phi(p_1),\cdots,\Phi(p_r))\in\mathbb P_{\Phi(\bomegas)}$ proving that $\Phi(\bpi)\in\wt_\ell V(\Phi(\bomega))$. For the converse let $(g_1,\cdots, g_r)\in\mathbb P_{\Phi(\bomegas)}$ be such that $\bomega(g_1)\cdots\bomega(g_r)\in\Phi(\cal P_{\bar i})$. Since the expression is reduced it follows that each $\bomega(g_s)\in\Phi(\cal P_{\bar i})$. By Lemma \ref{wtphi} we have  $(\Phi^{-1}(g_1),\cdots, \Phi^{-1}(g_r))\in\mathbb P_{\bomegas}$ and the converse follows.
\end{pf} 

 The next corollary is immediate using  Proposition \ref{weyl}.
\begin{cor}
For all $\bomega\in\cal P_{\bar i}^+$, we have  $$\bpi\in\wt_\ell W(\bomega)\implies \Phi(\bpi)\in\wt_\ell W(\Phi(\bomega)).$$
\end{cor}

  \subsubsection{The set $\wt_\ell V(\Phi(\bomega_{k,a}))\setminus\cal H_{\bar i, n}$} 
\begin{lem} \label{crucial} Suppose that  $k\in[1,\bar i]$ and $a\in\mathbb Z$ are such that $$\bomega_{ik, ia}\in\cal H_{\bar i,n},\ \ p\in\mathbb P_{ik,ia},\ \ \bomega(p)\in\wt_\ell V(\bomega_{ik,ia})\setminus\cal H_{\bar i, n}.$$
There exists $s\in\boc_p^-$ satisfying, $$\bomega_{s,p(s)}\notin \cal H_{\bar i, n}\ \ {\rm{and}}\ \ p(s)>p(j)\ \ {\rm{if}}\ \  j\in \boc_p^+\ \ {\rm{with}}\ \  \bomega_{j,p(j)}\notin\cal H_{\bar i, n}.$$   
\end{lem}
\begin{pf} Let $1\le j_1<\cdots<j_r\le n+1$ be an enumeration of $\boc_p^-\cup\boc_p^+$ and for $m\in[1,r]$ let $\epsilon_m=\pm 1 $ if $j_m\in\boc_p^\pm$. 
Write $\bomega(p)=\bomega_{j_1,p(j_1)}^{\epsilon_1}\cdots\bomega_{j_r, p(j_r)}^{\epsilon_r}$. Since $\bomega(p)\notin\cal H_{\bar i, n}$ it follows   that $$ r\ge 2,\ \ \boc_p^\pm\ne \emptyset,\ \ {\rm{and}}\ \  p(j_m)-p(j_{m-1})=\epsilon_m(j_{m-1}-j_{m}),\ \  m\in[2,r]. $$  
Suppose that $\bomega_{j_m,p(j_m)}\notin\cal H_{\bar i, n}$ and  $j_m\in\boc^+_p$. Then $j_{m-1}\in\boc_p^-$ if $m>1$  and $j_{m+1}\in\boc_p^-$ if $m<r$.  Setting $(j_0,\epsilon_0)= (0, 0)$ and   $(j_{r+1},\epsilon_{r+1})=(n+1,0)$ we have 
$$ p(j_{m-1})-p(j_m)=j_m-j_{m-1},\ \ p(j_{m+1})-p(j_m)=j_{m+1}-j_m.$$
If $\bomega_{j_{m\pm 1}, p(j_{m\pm 1})}\in\cal H_{\bar i,n}$ the first equation gives $p(j_m)+j_m=0\mod 2i$ while   the second equation gives  $p(j_m)-j_m=0\mod 2i$ which contradicts $\bomega_{j_m,p(j_m)}\notin\cal H_{\bar i, n}$. Hence at least of one of  $\bomega_{j_{m\pm 1}, p(j_{m\pm 1})}$ does not belong to $\cal H_{\bar i,n}$.
   Since $p(j_{m\pm 1})>p(j_m)$ it follows that $p(j_m)$ is not maximal  and the Lemma is proved.
\end{pf}

\subsubsection{} 
The next proposition gives a  partial converse to Corollary \ref{inflpaths}.
\begin{prop}\label{lweights2}  Suppose that $\bomega=\bomega_{j_1,a_1}\cdots\bomega_{j_k,a_k} \in\cal H_{\bar i,\bar i}^+$ and $
\bpi=\bomega(p_1)\cdots\bomega(p_k)\in\wt_\ell W(\Phi(\bomega))$ for some $ p_s\in\mathbb P_{ij_s,ia_s}$, $s\in[1,k]$. Then \begin{equation}\label{inflate1}\bpi\in\cal H_{\bar i,n}\cup\cal P_n^+\implies \bomega(p_s)=\Phi(\bomega(g_s)),\ \ g_s\in\mathbb P_{j_s,a_s},\  1\le s\le k.\end{equation}
\end{prop}
\begin{pf} In view of the Proposition \ref{inflpaths},  it suffices to prove that $\bomega(p_s)\in\cal H_{\bar i, n}$ for all $s\in[1,k]$. Suppose for   a contradiction that $\bomega(p_m)\notin\cal H_{\bar i,\bar n}$ and let $m$ be minimal with this property. Let  $s\in\boc_{p_m}^-$ be as in Lemma \ref{crucial}.  Since $\bpi\in\cal H_{\bar i, n}\cup\cal P_n^+$  there exists $m_1\in[1,k]$ with $s\in\boc_{p_{m_1}}^+$ and $p_{m_1}(s)=p_{m}(s)$. In particular this means that $\bomega(p_{m_1})\notin \cal H_{\bar i,n}$ and so $m_1>m$. Clearly this process can never stop which is  absurd
 and the proposition follows.
\end{pf}
\begin{cor}\label{ellwgtcor}
Suppose that $\bomega_s\in\cal H_{\bar i,\bar i}^+$ and $\bpi_s\in\wt_\ell(V(\Phi(\bomega_s))$ for $s=1,2$. Then $$\bpi_1\bpi_2\in\cal H_{\bar i,n}\cup\cal P_n^+\implies \bpi_s\in\cal H_{\bar i, n},\ \ s=1,2.$$
 \end{cor}
\begin{pf} 
Recall from Section \ref{weyl} that $V(\bomega)$ is the unique irreducible quotient of $W(\bomega)$ and that $$\wt_\ell V(\bomega)\subset \wt_\ell W(\bomega),\ \ \wt_\ell W(\bomega_1)\wt_\ell W(\bomega_2)=\wt_\ell W(\bomega_1\bomega_2).$$ 
Hence $\bpi_1\bpi_2\in\wt_\ell W(\Phi(\bomega_1)\Phi(\bomega_2))$ and the corollary is  now an immediate consequence of the proposition.
\end{pf}

\subsection{Proof of Theorem \ref{infliso}(i)--(iii)}

To prove that $\mathscr C_{\bar i,n}$ is a monoidal tensor category it suffices to prove that the Jordan--Holder constituents of $V(\bomega_1)\otimes V(\bomega_2)$ with $\bomega_s\in\cal H_{\bar i,n}^+$, $s=1,2$ are of the form  $V(\bpi)$ with $\bpi\in \cal H_{\bar i,n}^+$.  Writing $\bpi=\bpi_1\bpi_2$ with $\bpi_s\in\wt_\ell V(\bomega_s)$, $s=1,2$  it follows from  Corollary \ref{ellwgtcor} that $\bpi_s\in\cal H_{\bar i,n}$ for $s=1,2$ and hence $\bpi\in \cal H^+_{\bar i,n}$.\\\\
Let $\cal A$ be the subalgebra of   $\cal K_0(\mathscr C_{\bar i,n})$  generated by elements of the form $[V(\bomega_{ij,ia})]$ with $\bomega_{ij,ia}\in\cal H_{\bar i,n}$. Since $\cal K_0(\mathscr F_{n,\mathbb Z})$ is a polynomial algebra in the (infinitely many)  algebraically independent generators $[V(\bomega_{j,a})]$, $\bomega_{j,a}\in\cal P^+_{n}$ (see Theorem \ref{fr}) it follows that $\cal A$ is a polynomial algebra.\\\\
To prove that $\cal A= \cal K_0(\mathscr C_{\bar i,n})$ we show that $$\bomega\in\cal H_{\bar i ,n}\implies [V(\bomega)]\in\cal A$$  by using induction on the partial order on $P$, $\lambda\le\mu\iff\mu-\lambda\in Q_n^+$. The minimal dominant elements of the partial order are $0$ and $\{\omega_j: 1\le j\le n\}$ and induction
 obviously begins by definition of $\cal A$.
 For the inductive step we observe that Proposition \ref{weyl} gives  $$\bomega\in\cal H_{\bar i,n}\implies[ W(\bomega)]\in\cal A.$$
 On the other hand we can write $$[W(\bomega)]=\sum_{\bomegas'\in\cal P^+_n}a_{\bomegas'} [V(\bomega')]\in\cal K_0(\mathscr F_{n,\mathbb Z}),\ \ a_{\bomegas'}\in\mathbb Z.$$ Noting that $W(\bomega)_\mu\ne 0$ only if $\mu\le\wt\bomega$ and that $\dim W(\bomega)_{\wt\bomega}=1$ and 
using Corollary \ref{ellwgtcor} get $$ a_\bomegas=1,\ \   a_{\bomegas'}\ne 0\implies \bomega'\in\cal H_{\bar i, n}\ \  {\rm{and}}\ \  \bomega'\ne\bomega\implies\wt\bomega'<\wt\bomega.$$
  Hence the inductive hypothesis applies to all the  terms  with $\bomega'\ne\bomega$ and $a_{\bomegas'}\ne 0$. It follows that $[V(\bomega)]\in\cal A$ and part (i) of the theorem is proved. 
\\\\
For part (ii), note that  Corollary \ref{ellwgtcor} gives 
$$\bomega\in\cal H_{\bar i,n},\ \ \dim (V(\bomega_1)\otimes V(\bomega_2))_{\bomegas}\ne 0\implies \bomega=\bpi_1\bpi_2,\ \ \bpi_s\in\cal H_{\bar i,n}\cap\wt_\ell V(\bomega_s),\  s=1,2.$$ Together with the fact that $\chi^{\cal P_n}$ is  a homomorphism it follows that  for $\bomega_1,\bomega_2\in\cal H_{\bar i, n}$ and $\bomega\in\cal H_{\bar i, n}$ we have
$$\dim \left(V(\bomega_1)\otimes V(\bomega_2)\right)_\bomegas= \sum_{\bpis_1\bpis_2=\bomegas}\dim \left(V(\bomega_1)_{\bpis_1} \otimes V(\bomega_2)_{\bpis_2}\right),$$
proving that $\chi^{\cal H_{\bar i, n}}$ is a homomorphism.  
To prove that it is injective, 
recall from \cite[Lemma 6.17]{FM01}  that if $V_1$ and $V_2$ are two objects of $\mathscr F_{n,\mathbb Z}$ such that $\dim (V_1)_\bpis=\dim (V_2)_\bpis$, for all $\bpi\in \cal P_n^+$, then  $[V_1]=[V_2]$ in $\cal K_0(\mathscr F_{n,\mathbb Z})$. Since  Corollary \ref{ellwgtcor} implies that $\ell$--dominant weight of an object of $\mathscr C_{\bar i,n}$  is in $\cal H^+_{\bar i,n}$
  the injectivity is now clear.
The final statement of (ii) is immediate from \eqref{inflate1}.\\\\
For part (iii) the existence of the isomorphism $\Psi_{\bar i, n}$ sending $[V(\bomega_{j,a})]\mapsto [V(\bomega_{ij,ia})]$, $j\in[1,\bar i]$, $a\in\mathbb Z$ is clear since both algebras are polynomial on these algebraically independent generators. Using part (ii) of the theorem for $V$ an object of $\mathscr C_{\bar i,\bar i}$ we have, $$\chi^{\cal H_{\bar i, n}}(\Psi([V]))=\Phi\circ\ \chi^{\cal H_{\bar i,\bar i}}([V])=\sum_{\bomegas\in\cal H_{\bar i,\bar i}}\dim V_{\bomegas}e(\Phi(\bomega)).$$
Suppose that $\bomega\in\cal P_{\bar i}^+$ is a prime snake. By Proposition \ref{mysnake} we have $$\dim V(\bomega)_{\bomegas'}\in\{ 0,1\},\ \  \bomega'\in\cal P_{\bar i}. $$ By Proposition \ref{inflpaths} we have $$\dim V(\bomega)_{\bomegas'}= 1
\iff \dim V(\Phi(\bomega))_{\Phi(\bomegas')}= 1.$$ Hence, \begin{gather*}\chi^{\cal H_{\bar i,n}}([ V(\Phi(\bomega))])=\sum_{\bpis\in\cal H_{\bar i,n}}\dim V(\Phi(\bomega))_{\bpis} e(\bpi)= \sum_{\bpis\in\cal H_{\bar i,\bar i}}\dim V(\bomega)_{\bpis} e(\Phi(\bpi))\\ = \chi^{\cal H_{\bar i, n}}(\Psi([V(\bomega)])).\end{gather*}
Since $\chi^{\cal H_{\bar i,n}}$ is injective it follows that  $\Psi([V(\bomega)])=[V(\bomega)]$. If $\bomega$ is an arbitrary snake, then by \cite{MY12a} we can write $V(\bomega)$ as a tensor product of prime snakes. If  $\boa\in\mathbb Z^r$ is arbitrary with $\bomega_{j,\boa}\in\cal H_{\bar i, \bar i}$, then   we use Theorem \ref{main} to write $$V(\bomega_{j,\boa})\cong V(\bomega_{j,\boa_1})\otimes \cdots \otimes V(\bomega_{j,\boa_k}),$$ where $\boa_1,\cdots,\boa_k$ are $(j,\bar i)$--segments in general position. Clearly $\bomega_{j,\boa_s}\in\cal H_{\bar i, \bar i}$ for all $s\in[1,k]$. Using Proposition \ref{genpos} we see that the  $(ij,n)$--segments $i\boa_1,\cdots, i\boa_k$ are also in general position (since $a_j-a_{j-1}\in 2\mathbb Z$ by our hypothesis).
 Since $\Psi$ is a ring homomorphism we get 
 $$\Psi([V(\bomega_{j,\boa})])= \Psi([V(\bomega_{j,\boa_1})]) \cdots  \Psi([V(\bomega_{j,\boa_k})])=[V(\bomega_{ij,i\boa_1})] \cdots [V(\bomega_{ij,i\boa_k})],$$ and part (iii) is proved. 
 
\subsection{Monoidal Categorification and Proof of Theorem \ref{infliso}(iv)}  

The notion of a monoidal categorification of a cluster algebra was introduced in \cite{HL10} by Hernandez and Leclerc. They showed that the cluster algebra of type $A_{\bar i}$ is isomorphic to $\cal K_0(\mathscr C_{\bar i,\bar i}^1)$  via an isomorphism which maps a cluster monomial  to the class of an irreducible representation. \\\\
In terms of representation theory (see \cite{BC19a}) this amounts to solving the following problems in $\cal  K_0(\mathscr C_{\bar i,\bar i}^1)$:  classify the prime representations in the category $ \mathscr C_{\bar i,\bar i}^1$, give a necessary and sufficient condition for a tensor product of prime representations to be irreducible and show that the exchange relations hold for a pair of irreducible representations whose tensor product is reducible.  In  \cite{BC19a} we also describe the Jordan--Holder series of a reducible tensor product of prime representation.
Hence the equality in  \eqref{primeinf} in Theorem \ref{infliso}(iv) implies that that we also know the Jordan--Holder series of the image under $\Psi$. It is then immediate that $\mathscr C^1_{\bar i,n}$ is a monoidal tensor category and hence also that    $\cal  K_0(\mathscr C_{\bar i,n}^1)$ is  a monoidal categorification of type $A_{\bar i}$.\\\\
From now on our focus is to prove that $\Psi([V(\bomega)])=[V(\Phi(\bomega))]$ for all $\bomega\in\cal H_{\bar i,\bar i}^1$ and the proof occupies the rest of the section.
\subsubsection{} Set $I_{\bar i,n}= \{i,2i,\cdots, i\bar i\}\subset[1,n]$ and let 
$$
\epsilon_j=0,\ \ j\in 2i\mathbb Z,\ \ \epsilon_j=-i,\ \ 
 j\in  i(2\mathbb Z-1).$$ For $j,k\in I_{\bar i, n}$ define the following elements of $\cal H_{\bar i,n}$,$$\bof_j=\bomega_{j,3\epsilon_j}\bomega_{j, -\epsilon_j-2i},
  \ \ \bomega(j,k)=\bomega_{j,3\epsilon_j}\bomega_{j+i,3\epsilon_{j+i}}\cdots\bomega_{k,3\epsilon_k},\ \ j\le k,\ \ \bomega(j,k)=\bold 1,\ \ j> k.$$ 
  Set $$\bp\bor_{\bar i, n}=\{ \bof_j: j\in I_{\bar i,n}\}\sqcup \{ \bomega(j,k)  : j, k\in I_{\bar i,n},\ \ j\le k \}\sqcup\{\bomega_{m,-2i-\epsilon_m}: m\in I_{\bar i,n}\}.$$
Clearly
$$\bp\bor_{\bar i, n}=\Phi(\bp\bor_{\bar i, \bar i}).$$
 Given $\bomega=\bomega_{j_1,a_1}\cdots\bomega_{j_k,a_k}\in\cal P_n^+$  we set $\Ht\bomega=k$.
 \subsubsection{} The next proposition was proved in the special case when $i=1$ in \cite{HL10} (see also \cite{BC19a}). 
 \begin{prop} \label{irredcrit} Let $\bpi_1,\bpi_2\in\bp\bor_{\bar i,n}$. Then $$[V(\bpi_1)][V(\bpi_2)]=[V(\bpi_1\bpi_2)]$$ if one of the following hold:
 \begin{enumerit}
 \item[(i)] $\bpi_1=\bof_j$ for some $j\in I_{\bar i,n}$,
 \item[(ii)] $\bpi_1=\bomega(j,k)$ and $\bpi_2=\bomega(m,r)$ for $j,k,m,r\in I_{\bar i,n}$ with $j\le k$, $m\le r$ and $m\ne k+i$ and    one of the following hold:
 $$j=m,\ \ \ k\le m,\ \ k=r,\ \ \ j<m<k<r \ \ {\rm{and}} \ \epsilon_{k-m}=0,\ \ j<m<r<k\ \ {\rm{and}}\ \ \epsilon_{r-m}=-i.$$
 \item[(iii)] $\bpi_1=\bomega_{j,-\epsilon_j-2i}$ and $\bpi_2=\bomega_{k,-\epsilon_k-2i}$ or $\bpi_2= \bomega(m,r)$ with $j\notin[m,r]$.\end{enumerit}
\end{prop}
\begin{pf} 
 The proof of the proposition follows the lines of the proof of the corresponding result for $\mathscr C_{\bar i, \bar i}^1$ given in \cite[Theorem 3]{BC19a}. The basic idea is to proceed by induction on $\Ht\bpi_1\bpi_2$. We shall prove that induction begins; the proof of the inductive step is identical (and easier) to the one given in \cite[Section 3.7]{BC19a}.\\\\
To see that induction begins for part (i) we must prove that $V:= V(\bomega_{k,a})\otimes V(\bof_j)$ is irreducible for all $j,k\in I_{\bar i,n}$ and $a\in \{3\epsilon_k, -2i-\epsilon_k\}$. 
If $k=j$ the result is a particular case of Theorem \ref{main}. If $|j-k|>i$  it is easy to check that for $b\in \{3\epsilon_j, -2i-\epsilon_j\}$ we have
\begin{equation}\label{condirred}\pm (a-b) \notin \{2p+2-j-k: \max\{j,k\}\le p<\min\{j+k-1,n\}\}.\end{equation}
It follows from \cite{Ch01} that $V$ is irreducible.   If $|j-k|=i$ there are four cases to consider corresponding to  $\{\epsilon_k,\epsilon_j\}= \{0,-i\}$ and $a\in \{3\epsilon_k, -2i-\epsilon_k\}$. Assume first that $\epsilon_k =0$. If $a=-2i$ then  \eqref{condirred} again holds and the irreducibility  follows from \cite{Ch01} as before. If $a=0$ then  it suffices to prove that, 
$$\wt_\ell (V(\bomega_{k,0})\otimes V(\bof_j))\cap \cal P_n^+ =\{\bomega_{k,0}\bof_j\}.$$
To prove this we assume that $\bomega$ is an element of the intersection. Proposition \ref{inflpaths} and Corollary \ref{ellwgtcor} imply that we can write $\bomega=\Phi(\bomega_1)\Phi(\bomega_2)$ with $\bomega_1,\bomega_2\in \cal P_{\bar i}$,  $\bomega_1\in\wt_\ell V(\bomega_{k/i, 0}) $ and $\bomega_2\in V(\bof_{j/i})$.  Hence $\bomega_1\bomega_2\in\cal P^+_{\bar i }\cap \wt_\ell  V(\bomega_{k/i,0})\otimes V(\bof_{j/i})$ and it  follows from \cite[Proposition 4.3(a)]{HL13} that $\bomega_1\bomega_2=\bomega_{k/i,0}\bomega_{j/i,-1}\bomega_{j/i,-3}$.\\\\
 Suppose that $\epsilon_k=-i$. Then 
$$\Omega(V(\bomega_{k,a})\otimes V(\bof_j))^* = V(\bomega_{k,-a+n+1})\otimes V(\bomega_{j,-3\epsilon_j+n+1}\bomega_{j,2i+\epsilon_j+n+1}).$$
It follows from the previous case, after a suitable shift of parameters, that the module on the right is irreducible and hence, so is $V(\bomega_{k,a})\otimes V(\bof_j)$.\\\\

\noindent For parts (ii)-(iv)  to see that induction begins we have to establish the cases when $\Ht\bpi_1+\Ht\bpi_2\le 4$. Suppose that $\Ht\bpi_1+\Ht\bpi_2=2$. In this case, if we write
 $\bpi_1\bpi_2=\bomega_{j,a}\bomega_{k,b}$ it follows that $\bpi_1\bpi_2\notin \bp\bor_{\bar i,n}$.  Then $\pm (a-b)\notin \{2p+2-j-k: \max\{j,k\}\le p<\min\{j+k-1,n\}\}$ and it follows from \cite{Ch01} that $V(\bomega_{j,a})\otimes V(\bomega_{k,b})$ is irreducible.\\

Suppose that $\Ht\bpi_1+\Ht\bpi_2=3$. In this case  we have to prove that if  $\bomega(j,j+i)\bomega_{k,b}\notin \bp\bor_{\bar i,n}\cup\bof_k\bp\bor_{\bar i,n}$  then $V(\bomega(j,j+i))\otimes V(\bomega_{k,b})$ is irreducible. Note that the hypothesis implies that $$\bomega_{j,3\epsilon_j}\bomega_{k,b}\notin \bp\bor_{\bar i,n}, \ \ {\rm and} \ \ \  \bomega_{j+i,3\epsilon_{j+i}}\bomega_{k,b}\notin \bp\bor_{\bar i,n}.$$ If $\epsilon_{j}=-i$ and hence $\epsilon_{j+i}=0$ it follows from Proposition \ref{weyl}(iii) and the preceding case that 
$$V(\bomega_{j+i,0})\otimes V(\bomega_{j,-3i})\otimes V(\bomega_{k,b}), \ \ V(\bomega_{k,b})\otimes V(\bomega_{j+i,0})\otimes V(\bomega_{j,-3i})$$
are both $\ell$-highest weight modules. Hence $V(\bomega(j,j+i))\otimes V(\bomega_{k,b})$ and its dual 
are $\ell$--highest weight and thus irreducible by Proposition \ref{weyl}(i).
 The case $\epsilon_{j}=0$ is completely analogous and  we omit the details.\\
 
If $\Ht\bpi_1+\Ht\bpi_2=4$ then we have to prove that the following modules are irreducible:
$$V(\bomega(j,j+i))\otimes V(\bomega(j,j+i)), \ \ 
 V(\bomega(j,j+i))\otimes V(\bomega(j+i,j+2i)).$$
 The first assertion is immediate from \cite{DLL19} since $\bomega(j,j+i)$ is a prime snake. For the second, 
 if $\epsilon_j =0$ then 
 $$\bomega(j,j+i)=\bomega_{j,0}\bomega_{j+i,-3i}, \ \ \bomega(j+i,j+2i)=\bomega_{j+i,-3i}\bomega_{j+2i,0}.$$
Using Proposition \ref{weyl}(iii) and the previous case  it follows that the modules 
$$V(\bomega_{j,0}) \otimes V(\bomega_{j+i,-3i}) \otimes V(\bomega(j+i,j+2i)), \ \  V(\bomega_{j+2i,0})\otimes V(\bomega_{j+i,-3i}) \otimes V(\bomega(j,j+i))$$
are $\ell$-highest weight and the irreducibility of   $V(\bomega(j,j+i))\otimes V(\bomega(j+i,j+2i))$ follows as before. The proof when $\epsilon_j=-i$ is entirely analogous so we omit details.
\end{pf}

 \subsubsection{}
 \begin{prop} \label{indstep} For  $p, j,k\in I_{\bar i,n}$ with $j<k$, the following identities hold in $\cal K_0(\mathscr C_{\bar i,n}^1)$.
$$(i)\ \ [V(\bomega_{p,3\epsilon_p})][V(\bomega_{p,-2i-\epsilon_p})]=[V(\bof_p)]+[V(\bomega_{p-i,\epsilon_p-i})][V(\bomega_{p+i,\epsilon_p-i})].$$
  $$(ii)\ \ [V(\bomega_{j,-2i-\epsilon_j})][V(\bomega(j,k))] = [V(\bof_j)][V(\bomega(j+i,k))] + [V(\bomega_{j-i,\epsilon_j-i})][V(\bof_{j+i}\bomega(j+2i,k))],$$ 
$$(iii)\ \ [V(\bomega_{j,3\epsilon_j})][ V(\bomega(j+i,k))]=  [V(\bomega(j,k))]+ [V(\bpi)],$$ where $$[V(\bpi)]=[ V(\bomega_{j-i,\epsilon_j-i})][ V(\bomega_{j+2i,\epsilon_{j+i}-i})]^{\delta_{j+i,k}}[V(\bof_{j+2i} \bomega(j+3i,k))]^{1-\delta_{j+i,k}},$$

 \end{prop}
\begin{pf}  For $j,k\in I_{\bar i,n}$ with $k\ge j$, set
\begin{gather*}M(j,k)=V(\bomega_{j,3\epsilon_j})\otimes V(\bomega(j+i,k)) \otimes V(\bomega_{j,-2i-\epsilon_j})\ \ \  \\ U(j,k)=V(\bomega(j,k))\otimes V(\bomega_{j,-2i-\epsilon_j}),\ \ 
W(j,k)= V(\bof_j)\otimes V(\bomega(j+i,k)),\end{gather*}
Then parts (i) and (ii) the proposition can be formulated as:
\begin{equation}\label{parti} [M(j,j)]=  [U(j,j)]= [W(j,j)]+ [V(\bomega_{j+i,\epsilon_j-i})][V(\bomega_{j-i,\epsilon_j-i})]
\end{equation}
\begin{equation}\label{partii}[U(j,k)] = [W(j,k)]+[ V(\bomega_{j-i,\epsilon_j-i})][W(j+i,k)],\ \ k>j.\end{equation} Since the Grothendieck ring has no zero divisors we see that the following equations are equivalent to part (iii) of the proposition:
\begin{equation}\label{partiiia}[M(j,j+i)]= [U(j,j+i)]+[V(\bomega_{j-i,\epsilon_{j}-i})][V(\bomega_{j+2i,\epsilon_{j+i}-i})][V(\bomega_{j, -2i-\epsilon_j})],\end{equation}
\begin{equation}\label{partiiib}[M(j,k)]= [U(j,k)]+[ V(\bomega_{j-i,\epsilon_j-i})] [V(\bomega_{j,-2i-\epsilon_j})][W(j+2i,k)],\ \ k\ge j+2i. \end{equation}
\vskip12pt
\noindent The equalities in \eqref{parti} and \eqref{partiiia} are immediate from Proposition \ref{mysnake}.  Since  the module $V(\bomega(j,k))$ occurs in the Jordan--Holder series of the reducible tensor product   $V(\bomega(j+i,k))\otimes V(\bomega_{j,3\epsilon_j})$  we can write   $$ [M(j,k)]= [U(j,k)]+[K_1(j,k)], \ \ k\ge j$$ where $[K_1(j,j)]=0$ and  $[K_1(j,k)]$ is the class of an object of $\mathscr F_{n,\mathbb Z}$. Note that by \eqref{partiiia} we have $$[K_1(j,j+i)]=[V(\bomega_{j-i,\epsilon_{j}-i})][V(\bomega_{j+2i,\epsilon_{j+i}-i})][V(\bomega_{j, -2i-\epsilon_j})].$$
The module $W(j,k)$ is irreducible by Proposition \ref{irredcrit} and hence occurs in the Jordan--Holder series of both $M(j,k)$ and $U(j,k)$. In particular we can write $$[M(j,k)]=[W(j,k)]+[K_2(j,k)],\ \ [U(j,k)]=[W(j,k)]+[K_3(j,k)],\ \  k\geq j$$ where $K_s(j,k)$, $j=2,3$ is the class of an object of $\mathscr F_{n,\mathbb Z}$.
Combining the  expressions for $[M(j,k)]$ and $[U(j,k)]$ we get $$[K_1(j,k)]+ [K_3(j,k)]= [K_2(j,k)].$$ Note that $$[K_2(j,j)]= [K_3(j,j)]=[V(\bomega_{j-i,\epsilon_j-i})][V(\bomega_{j+i,\epsilon_j-i})].$$
Proposition \ref{mysnake} gives \begin{equation}\label{k2} [K_2(j,k)]= [V(\bomega_{j-i,\;\epsilon_j-i})][ V(\bomega_{j+i,\;\epsilon_j-i})][ V(\bomega(j+i,k))]=[V(\bomega_{j-i,\epsilon_j-i})][ U(j+i,k)],\end{equation} where the last equality holds if $k>j$.
Hence we have 
 \begin{equation}\label{crux1}
 [K_3(j,j+i)]+[K_1(j,j+i)]= [V(\bomega_{j-i,\epsilon_j-i})]\left([V(\bof_{j+i})]+ [V(\bomega_{j,-2i-\epsilon_j})][V(\bomega_{j+2i,\epsilon_{j+i}-i})]\right),    
 \end{equation}
\begin{equation}\label{crux}[K_3(j,k)]+[K_1(j,k)]=[V(\bomega_{j-i,\epsilon_j-i})]\left([K_3(j+i,k)]+   [W(j+i,k)]\right) ,\ \ k\ge j+2i.\end{equation} 
The equalities in  \eqref{partii} and \eqref{partiiib} hold if we prove that $$[K_1(j,k)]= [V(\bomega_{j-i,\epsilon_j-i})][ K_3(j+i,k)],$$ $$[K_3(j,j+i)]= [V(\bomega_{j,\epsilon_j-i})][\bof_{j+i}],\ \  [K_3(j,k)]= [V(\bomega_{j-i,\epsilon_j-i})][ W(j+i,k)],\ \ k\ge j+2i.$$ We proceed by 
 an induction on $k$. To see that induction begins at $k=j+i$, we 
observe that  \eqref{crux1} and  \eqref{partiiia} give $[K_3(j,j+i)]= [V(\bomega_{j-i,\epsilon_j-i})][V(\bof_{j+i})]$ as needed. \\\\
 Suppose that $k\ge j+2i$.  Note that the $\ell$--weight $\bomega_{j-i,\epsilon_j-i}\bof_{j+i}\bomega(j+2i,k)$ of $V(\bomega_{j-i,\epsilon_j-i})\otimes W(j+i,k)$ occurs in both $M(j,k)$ and $U(j,k)$ with multiplicity one so it cannot occur in $K_1(j,k)$. By Proposition \ref{irredcrit} each of the terms on the right hand side of \eqref{crux} is the class of an irreducible representation. It follows that
$$[K_1(j,k)]= [V(\bomega_{j-i,\epsilon_j-i})][ K_3(j+i,k)],\ \ [K_3(j,k)]= [V(\bomega_{j-i,\epsilon_j-i})][ W(j+i,k)].$$ Using the inductive hypothesis we get the desired formulae.\end{pf}

 \subsubsection{Proof of Theorem \ref{infliso}(iv)} Notice that by Theorem \ref{infliso}(iii) we know part (iv) for the elements $\bomega\in \bp\bor_{\bar i,\bar i}$  such that $\Ht\bomega\leq 2$. 
A straightforward induction on $\Ht\bomega(j,k)$ along with Proposition \ref{indstep} (iii) show that part (iv) holds for all $\bomega(j,k)\in \bp\bor_{\bar i,\bar i}$.\\\\
To complete the proof of part (iv) it suffices to show that if $\bomega,\bomega'\in\bp\bor_{\bar i,\bar i}$ are such that $V(\bomega)\otimes V(\bomega')$ is irreducible then so is
$V(\Phi(\bomega))\otimes V(\Phi(\bomega')).$ For this we recall that in \cite{BC19a} we actually proved that in the category $\mathscr C_{\bar i,\bar i}^1$, the conditions for irreducibility given in Proposition \ref{irredcrit} are also necessary. Observing that if $\bomega,\bomega'\in\bp\bor_{\bar i,\bar i}$ satisfy the conditions in Proposition \ref{irredcrit} then  $\Phi(\bomega),\Phi(\bomega')\in\bp\bor_{\bar i,n}$ satisfy the conditions in Proposition \ref{irredcrit} and hence $V(\Phi(\bomega))\otimes V(\Phi(\bomega'))$ is irreducible and the proof is complete.\hfill\qedsymbol

\begin{rem}  In \cite[Theorem 7.8]{HL10} it is proved that for each $\bomega\in \bp\bor_{\bar i,\bar i}$ we have 
$$\chi^{\cal H_{\bar i,\bar i}}V(\bomega) = \bomega F_\bomegas(\balpha_{1,\epsilon_1-1}^{-1},\cdots, \balpha_{\bar i,\epsilon_{\bar i}-1}^{-1}),$$
where $F_\bomegas(y_1,\cdots, y_{\bar i})$ is the $F$--polynomial of the cluster variable associated to $[V(\bomega)]$.  Therefore, using  \eqref{phichi} we have that 
$$\chi^{\cal H_{\bar i,n}}V(\Phi(\bomega)) = \Phi\circ \chi^{\cal H_{\bar i,\bar i}}V(\bomega) = \Phi(\bomega) F_\bomegas (\Phi(\balpha_{1,\epsilon_1-1}^{-1}),\cdots, \Phi(\balpha_{\bar i,\epsilon_{\bar i}-1}^{-1})),$$
(see Remark \ref{infroot}). In this sense, $F_\bomegas$ can be identified with the evaluation of the $F$--polynomial of $[V(\Phi(\bomega))]$, however the evaluation is different from the one in \cite{HL10}.
\end{rem}

\section{Imaginary Modules}\label{imaginary}
 Recall that an irreducible  $\hat\bu_n$--module  is said to be real if its  tensor square is irreducible and  otherwise the module is said  to be  imaginary. It has been known since 1991 (see \cite{CP91}) that imaginary modules do not exist when $n=1$. The first example of an   imaginary module was  given in 2006 by Leclerc in \cite{Lec02}, where he proved that if  $n=5$,  the module $V(\bomega)$  with  $\bomega=\bomega_{2,6}\bomega_{1,3}\bomega_{3,3}\bomega_{2,0}\in\cal P^+_5$  is imaginary. Subsequently further  examples of  imaginary modules  were given  in \cite{LM18}.  In both cases, the proof that the modules are imaginary is indirect; examples are first constructed for the  affine Hecke algebra and then affine Schur--Weyl duality gives that the corresponding module for  a quantum affine algebra $\hat{\bu}_N$ for some $N$ sufficiently large is imaginary.  Note that the reducibility of a tensor product in rank $N$ does not necessarily  imply reducibility for all $r<N$. For instance,  as we shall see in this section,  the example given by Leclerc, is actually imaginary for $\hat\bu_3$ but this cannot be deduced by using  Schur--Weyl duality if $n\le 3$.  It was noted in \cite[Section 13.6]{HL10} that the example given by Leclerc, is actually imaginary for $\hat\bu_3$ and in the language of the current paper is  an object of $\mathscr C_{3,3}^3$.  But no examples of imaginary modules were known in types $A_2$ and we  prove the existence of imaginary modules  in this case as well.
Here, we remind the reader that Hernandez and Leclerc speculated imaginary  simple objects existed in $\mathscr C_{n,n}^\ell$ whenever $\cal A_\ell$ does not have a finite cluster type.  In the case of type $A_2$ this means that the category $\mathscr C_{2,2}^5$ (c.f. \cite[Table 1]{HL10}) should have imaginary modules and our smallest example of imaginary modules in this case is exactly  an object of this category.  For $n\ge 2$  our examples live in $\mathscr C_{n,n}^\ell$ for $\ell\ge n$.\\\\
  In this section we shall give a systematic way to construct imaginary modules by working with tensor products of $\hat\bu_n$--modules associated to $(i,n)$--segments with $n\ge 2i$ and $i\ge 2$. In Section \ref{a2} we  give  a family of examples of  imaginary modules for the quantum affine algebra associated to $A_2$. In Section \ref{d4}, we construct a rather the first examples of  imaginary modules in $D_4$ which do not arise from an $A_n$--example. In fact we shall see that the corresponding $A_3$--module is real.

 \subsection{} The following is  the main result of this section.
\begin{thm}\label{imagen} Let $n\ge 3$ and $i\in[2,n-1]$ with $2i\le n+1$ and $r\geq 2$.  Let $\bob = (b_1,\cdots, b_r)$ be an $(i,n)$--segment with $b_j-b_{j-1}<2i$ for all $j\in[2,r]$ and  let $$\boa=(a_1,\cdots, a_r),\ \ a_k= b_k-n-1,\ \ k\in[1,r],\ \  s_j=\frac12(b_{j-1}-b_j+2i),\ \ j\in[2,r].$$  Setting
 $$\bomega=\bomega_{i, b_r}(\bomega_{i-s_{r}, b_{r-1}-s_{r}}\bomega_{n+1-i+s_{r}, a_r+s_r})\cdots(\bomega_{i-s_2, b_1-s_2}\bomega_{n+1-i+s_2, a_2+s_2})\bomega_{n-i+1, a_1},$$ the module $V(\bomega)$
 is imaginary.
 \end{thm}
 The proof of the theorem is given in Sections \ref{lwtexists}--\ref{lwtpf}.
 
 \subsection{The key step}
 The proof of the following proposition can be found in Section \ref{lwtpf}. 
\begin{prop}\label{lwtexists} Retain the notation established so far and 
 let $\bpi\in\wt_\ell^+ W(\bomega)$. Then $$\Hom_{\hat\bu_n}(W(\bpi), V(\bomega_{i,\bob})\otimes V(\bomega_{n+1-i,\boa}))\ne 0\iff \bpi\in\{\bold 1, \bomega\}.$$ 
\end{prop}
\subsection{Proof of Theorem  \ref{imagen}}  Assuming Proposition \ref{lwtexists} we complete the proof of the theorem.

 \subsubsection{The map $\Phi$}\label{phi} Noting that $\bomega_{n+1-i,\boa}={}^*\bomega_{i,\bob}$ we have   $\hat\bu_n$--maps,
\begin{equation}\label{sochd}\mathbb C\hookrightarrow V(\bomega_{i,
 \bob})\otimes V(\bomega_{n+1-i,\boa}),\ \ <\ \ >: V(\bomega_{n+1-i,
 \boa})\otimes V(\bomega_{i,\bob})\to \mathbb C\to 0.\end{equation}
 Let $$\Phi:V(\bomega_{i,
 \bob})\otimes V(\bomega_{n+1-i,\boa})\otimes V(\bomega_{i,
 \bob})\otimes V(\bomega_{n+1-i,\boa})\to V(\bomega_{i,
 \bob})\otimes V(\bomega_{n+1-i,\boa})\to 0,$$ be the composition $\id\circ<\ \ >\circ\id.$\\\\ 
 Since   $V(\bomega_{i,\bob})$ is real by Theorem \ref{main} it follows that $\mathbb C$ is the socle (resp. head) of the tensor products in \eqref{sochd}. In particular any submodule  of $V(\bomega_{i,\bob})\otimes V(\bomega_{n+1-i,\boa})$ must contain the trivial representation. 
\noindent In what follows, we  shall make use of this remark without further mention.

 \subsubsection{}   We need the following result on the set of $\ell$--weights of a tensor product of a prime snake module and its dual. 
 \begin{lem}\label{dimv} Let $\bpi=\bomega_{i_1, c_1}\cdots\bomega_{i_k, c_k}$  with $i_1,\cdots, i_k\in[1,n]$ and $c_1<\cdots<c_k$  be a prime snake and set $V:= V(\bpi)\otimes V(^*\bpi)$. Then  $$\dim V_{\bold 1}= 1= \dim (V\otimes V)_{\bold 1}.$$ Moreover if $\bpi_1\in \cal P_n$ has a reduced expression containing $\bomega_{n+1-i_1,c_1-n-1}^{s}$ for some  $s\ge 2$ then $\bpi_1\notin\wt_\ell V.$

 \end{lem}
 \begin{pf}
 Suppose that $\tilde\bpi\tilde \bpi'=\bold 1$ for some $\tilde\bpi\in  \wt_\ell V(\bpi)$ and $\tilde \bpi'\in\wt_\ell V({}^*\bpi)$ and write $$\tilde \bpi'=\bomega(g_1')\cdots \bomega(g_k'),\ \ \tilde \bpi=\bomega(g_1)\cdots \bomega(g_k),\ \ (g_1',\cdots, g'_k)\in\mathbb P_{{}^*\bpis}, \ (g_1,\cdots, g_k)\in\mathbb P_{\bpis}.$$ Since $\wt_\ell V(\bpi)$ is in the subgroup of $\cal P$ generated by elements $\bomega_{j,c}$ with $c\ge c_1$ we see that $\tilde\bpi'$ must also be in this subgroup.
  It follows that $\bomega(g_1')=\bomega_{i_1,c_1}^{-1}$. Since $\bomega_{i_1,c_1}$ cannot occur in $\bomega(g_s)$ for any $s>1$ we get $\bomega(g_1)=\bomega_{i,c_1}$. An obvious iteration proves that $\tilde\bpi'=\bomega_{i,\bob}^{-1}$ and hence $\dim V_{\bold 1}=1$. The proof that $\dim (V\otimes V)_{\bold 1}=1 $ is identical. The second assertion is immediate from Proposition \ref{mysnake} applied to ${}^*\bpi$ since $\bomega_{n+1-i_1,c_1-n-1}$ cannot occur in any $\ell$--weight of $V(\bpi)$.

 \end{pf}

\subsubsection{} 
\begin{lem}\label{lenth2}
There exists an $\ell$--highest submodule $M(\bomega)$ of $V(\bomega_{i,\bob})\otimes V(\bomega_{n+1-i,\boa})$ and a non--split short exact sequence  $$0\to\mathbb C\to M(\bomega)\to V(\bomega)\to 0.$$
\end{lem}
\begin{pf} By Proposition \ref{lwtexists} there exists a non--zero map of $\hat\bu_n$--modules $$f:W(\bomega)\to  V(\bomega_{i,\bob})\otimes V(\bomega_{n+1-i,\boa}),\ \ f(v_\bomegas)\ne 0 .$$ Then $M(\bomega)=\hat\bu_n f(v_\bomegas)$ is an $\ell$--highest weight module and hence so is any quotient of it. Moreover,   $V(\bomega)$ is  its unique irreducible quotient and since  $\mathbb C\hookrightarrow M(\bomega)$ it follows that $M(\bomega)/\mathbb C$ is $\ell$--highest weight. Suppose that $M(\bomega)/\mathbb C$ is reducible  and let $\bar N$ be the maximal submodule of $M(\bomega)/\mathbb C$. Taking $N$ to be the preimage in $M(\bomega)$ of $\bar N$ we see that $\bar N$ must contain an $\ell$--highest weight vector with $\ell$-weight $\bpi\notin\{\bold 1, \bomega\}$.
 But this is impossible by the forward direction of Proposition \ref{lwtexists} and proves that $M(\bomega)/\mathbb C$ is irreducible. The fact that the short exact is non--split is immediate since $M(\bomega)$ is $\ell$--highest weight  and the  lemma is proved. 
\end{pf}
\subsubsection{}    The next proposition and its corollary  establishes Theorem \ref{imagen}. 
\begin{prop}
The restriction of $\Phi$ to $M(\bomega)\otimes M(\bomega)$ is non--zero and $\Phi(M(\bomega)\otimes M(\bomega))$ has an irreducible quotient $V(\bpi)$ with $\bpi\in\cal P_n^+\setminus\{\bold 1, \bomega,\bomega^2\}$.
\end{prop}
\begin{pf} Set $V:= V(\bomega_{i,\bob})\otimes V(\bomega_{n+1-i,\boa})$.
Using Lemma \ref{dimv} and Lemma \ref{lenth2} we have  $$\dim(M(\bomega)\otimes M(\bomega))_{\bold 1}=\dim (V\otimes V)_{\bold 1}=\dim V_{\bold 1} =1. $$ Since $\Phi: V\otimes V\to V$ is surjective it follows that $\Phi(M(\bomega)\otimes M(\bomega))\ne 0$. Since $\bomega^2\notin\wt_\ell V$ by Lemma \ref{dimv} it follows that $V(\bomega^2)$ is  not  a quotient of $\Phi(M(\bomega)\otimes M(\bomega))$. 
 We claim that $\mathbb C$ is also not a quotient of $M(\bomega)\otimes M(\bomega)$. Otherwise we would have a non--zero element, say $f\in\Hom_{\hat\bu_n}(M(\bomega), {}^*M(\bomega))$. By Lemma \ref{lenth2} we have $$0\to V(^{*}\bomega)\to {}^*M(\bomega)\to \mathbb C\to 0.$$ Since $M(\bomega)$ is $\ell$--highest weight it follows that the  image of $f$ must be isomorphic to $V(^{*}\bomega)$. But this is absurd since $\bomega\ne {^{*}}\bomega$. 
Hence  $\mathbb C$ is  also not a quotient of $\Phi(M(\bomega)\otimes M(\bomega))$. 
To prove that $V(\bomega)$ is not a quotient as well, it suffices to show (since $M(\bomega)$ is an $\ell$--highest weight module) that $$\Hom_{\hat\bu_n}(W(\bomega)\otimes W(\bomega), V(\bomega))=0.$$ 
This follows from  Lemma \ref{lowestwt}  if we prove that $(\bomega^*)^{-1}\bomega\notin\wt_\ell W(\bomega).$  Write, $$\bomega=\bomega_1\bomega_2,\ \ \bomega_1=\bomega_{n+1-i,a_1}\prod_{j=2}^m\bomega_{i-s_j,b_1-s_j}\prod_{j=2}^{p}\bomega_{n+1-i+s_j,a_j+s_j}, $$
where $m=\max\{j\in[2,r]: b_{j-1}-s_j<b_1\}$ and $p=\max\{j\in[2,r]: a_j+s_j<b_1\}$.  Notice that $\bomega_1\bomega_{i,b_1}^{-1}$
occurs in a reduced expression for $(\bomega^*)^{-1}\bomega$. We have $\wt_\ell W(\bomega)=\wt_\ell W(\bomega_1)\wt_\ell W(\bomega_2)$ and also that $\wt_\ell W(\bomega_2)$ is contained in the submonoid generated by $\bomega_{k,c}^{\pm 1}$ with $c>b_1$ and $\bomega_{k,b_1}$ for $k\in[1,n]$. This means that $\bomega_1\bomega_{i,b_1}^{-1}$
must be in $\wt_\ell W(\bomega_1)$ which is clearly impossible.
\end{pf}

\begin{cor}
There exists $\bpi\in\cal P_n^+$ with $\bpi\ne\bomega^2$ and a non--zero map $V(\bomega)\otimes V(\bomega)\to V(\bpi)$. In particular, $V(\bomega)$ is imaginary.
\end{cor}
\begin{pf} The proposition implies that there exists $$f:\Phi(M(\bomega)\otimes M(\bomega))\to V(\bpi)\to 0,\ \ \bpi\in\cal P_n^+\setminus\{\bold 1, \bomega,\bomega^2\}.$$
Noting that  $$V(\bomega)\otimes V(\bomega)\cong \left(M(\bomega)\otimes M(\bomega)\right)/(M(\bomega)\otimes \mathbb C+\mathbb C\otimes M(\bomega)), $$  $$\bpi\notin\{ \bold 1,\bomega\}\implies  f(\Phi(M(\bomega)\otimes\mathbb C+\mathbb C\otimes M(\bomega)))=0,$$  we have an induced map $V(\bomega)\otimes V(\bomega)\to V(\bpi)\to 0$. Since $\bpi\ne\bomega^2$ this map cannot be an isomorphism and the corollary follows.
\end{pf}

 \subsection{Proof of Proposition \ref{lwtexists}}\label{lwtpf}
 \subsubsection{} 
 A simple calculation shows that $$\bomega=\bomega_{i,\bob}\bomega(\underline p),\ \ \underline p=(g_{n+1-i,a_1}^0, g_{n+1-i,a_2}^{s_2},\cdots , g_{n+1-i,a_r}^{s_r})\in \mathbb P_{n+1-i,\boa}.$$
 Assume that $\Hom_{\hat{\bu}_n}(W(\bpi),V(\bomega_{i,
 \bob})\otimes V(\bomega_{n+1-i,\boa}))\ne 0 $ for some $\bpi\in\cal P_n^+$. Using   Lemma \ref{socle} and  Section \ref{gjm} we write$$\bpi=\bomega_{i,\bob}\bomega(\underline g),\ \  \underline g=(g_{n+1-i,a_1}^{m_1},\cdots, g_{n+1-i,a_r}^{m_r})\in\mathbb P_{n+1-i, \boa},\ \ m_j\in[0,i],\ \ j\in[1,r],$$
\begin{equation}\label{mj}\ \ {\rm{where}}\ \ m_j>0\implies b_j-2i+2m_j\in\{b_1,\cdots, b_j\}.\end{equation}
 If in addition $\bpi\in\wt_\ell W(\bomega)$ then $\bomega(\underline p)\in\bomega(\underline g)\cal Q^+,$ and an application of Proposition \ref{mysnake}(iv) now gives $m_j\ge s_j>0$ for all $j\in[2,r]$.
 Since
   $$i\ge m_j\ge s_j\implies b_{j-1}=b_j+2(s_j-i)\le b_j+2(m_j-i)\le b_j,$$ we see that \eqref{mj} forces $b_j-2i+2m_j=b_j$, i.e. $m_j\in\{s_j, i\}$ for all $j\in[2,r]$. If $j=1$ then again the restriction that $\bpi\in \cal P_n^+$ implies $m_1\in\{0,i\}$.
Finally note that these arguments prove that  $m_j=i$ implies $m_{j+1}=i$.
   Summarizing, we have proved that if $\bpi\notin\{\bold 1,\bomega\}$ then
\begin{equation}\label{bpiposs}\bpi=\bomega_{i,b_j}(\bomega_{i-s_j, b_{j-1}-s_j}\bomega_{n+1-i+s_j, a_j+s_j})\cdots (\bomega_{i-s_2,b_1-s_2}\bomega_{n+1-i+s_2,a_2+s_2})
 \bomega_{n+1-i,a_1}\end{equation}  for some $j\in[1,r-1]$. 
 The second condition in Lemma \ref{socle} gives $$(\bpi^*)^{-1}\bomega_{i,\bob}=\bomega(g_1)\cdots\bomega(g_r),\ \ (g_1,\cdots, g_r)\in\mathbb P_{i,\bob}.$$
 A reduced expression of the left hand side of the preceding equality always contains $\bomega_{i,b_k}$ for all $k\in[2,r]$ and does not contain $\bomega_{n+1-i, b_k +n+1}^{-1}$ for any $k\ne j$. Hence a standard argument using Proposition \ref{mysnake} shows that $\bomega(g_p)=\bomega_{i,b_p}$ for all $p\in [2,r]$ with $p\ne j$ and $\bomega(g_j)\in\{\bomega_{i,b_j}, \bomega_{i,b_j+n+1}^{-1}\}$.
 But this means that $$\bomega(g_1)\cdots\bomega(g_r)=\bomega_{i,b_r}\cdots\bomega_{i,b_j}^\epsilon\bomega_{n+1-i,b_j+n+1}^{-1+\epsilon}\cdots\bomega_{i,b_2}\bomega(g_1),\ \ \epsilon=0,1$$ or equivalently that $$(\bpi^*)^{-1}\bomega_{i,b_j}\bomega_{i,b_1}= \bomega_{i,b_j}^\epsilon\bomega_{n+1-i, b_j+n+1}^{-1+\epsilon} \bomega(g_1).$$ Substituting for $\bpi$ we are forced to have $\epsilon=0$. i.e. $\bomega(g_j)=\bomega_{n+1-i, b_j+n+1}^{-1}$. But this means that $g_j(n+1-i)=b_j+n+1> g_{j+1}(n+1-i)$ and hence we have a contradiction to $(g_1,\cdots, g_r)\in\mathbb P_{i,\bob}$.

 \subsubsection{} The following proposition is needed for the converse direction. 
 \begin{prop} \label{dsub} Let  $J=\{j+1,\cdots ,k-1\}$ be a subset of $I$ and assume that
   $\bpi$ (resp. $\bpi_1,\bpi_2$) are in the submonoid of $\cal P_n^+$ generated by $\bomega_{j,c}$ with $j\in J$ (resp.  $I\setminus J$).
Setting \begin{gather*} \bpi=\bomega_{j_1,c_1}\cdots\bomega_{j_s,c_s},\\ \bpi' =\bomega_{k+j-j_1,c_1+k-j}\cdots\bomega_{k+j-j_s,c_s+k-j}, \ \ \ '\bpi =\bomega_{k+j-j_1,c_1-k+j}\cdots\bomega_{k+j-j_s,c_s-k+j}.\end{gather*}
 we have non--zero maps \begin{gather*}W(\bpi_1\bpi^+\bpi_2)\to  V(\bpi_1\bpi')\otimes V(\bpi\bpi_2),\ \ W(\bpi_1\bpi^-\bpi_2)\to  V(\bpi_1\bpi)\otimes V('\bpi\bpi_2),\\ {\rm{where,}}\ \  \bpi^\pm =(\bomega_{j,c_1\pm(j_1-j)}\bomega_{k,c_1\pm(k-j_1)})\cdots (\bomega_{j,c_s\pm(j_s-j)}\bomega_{k,c_s\pm(k-j_s)}).\end{gather*}
 \end{prop}
 \begin{pf}
 Let $v_1\in V(\bpi_1\bpi')$ and $v_2\in V(\bpi\bpi_2)$ be the highest $\ell$-weight vectors. By Section \ref{diagsub} we have the following  maps of $\hat\bu_{n,J}$--modules
 $$V_J(\bpi)\hookrightarrow V(\bpi_s\bpi),\ \ s=1,2,\ \  \hat\bu_{n,J}\ v_2 \cong V_J(\bpi),\ \ \ \ V_J(\bpi)^*\cong \hat\bu_{n,J}\ v_1 \cong V_J(\bpi').$$ Hence we  have non--zero maps of $\hat\bu_{n,J}$--modules,$$\mathbb C\hookrightarrow V_J(\bpi')\otimes V_J(\bpi)\hookrightarrow V(\bpi_1\bpi')\otimes V(\bpi_2\bpi),$$   and we let $v\in V(\bpi_1\bpi')\otimes V(\bpi_2\bpi)$ be the image of $1\in \mathbb C$.
 It is straightforward to compute that,
 $$\bpi'\bpi=\balpha_J,\ \ \balpha= \prod_{p=1}^s\prod_{r=j_p}^{k-1}\prod_{t=0}^{j_p-j-1}\balpha_{r-t,c_p+t+1}.$$
 Section \ref{diagsub} implies that $v$ is an $\ell$--highest weight vector with $\ell$--weight $\bpi_1\bpi_2\bpi'\bpi\balpha^{-1}$  and a further simple computation shows that $\bpi'\bpi\balpha^{-1}=\bpi^+$ which proves the existence of the first map.
The existence of the second map in the statement of the proposition is proved in an identical fashion and we omit the details. \end{pf}
 
 \subsubsection{} Set, $$\bomega_0=\bomega_{n+1-s_2, a_2-i+s_2}\cdots \bomega_{n+1-s_r, a_r-i+s_r},$$ and write $\bomega=\bomega_1\bomega_2$, where 
 $$\bomega_1=\bomega_{i,b_r}\bomega_{i-s_r, b_{r-1}-s_r}\cdots\bomega_{i-s_2, b_1-s_2},\ \ \bomega_2=\bomega_{n+1-i+s_r, a_r+s_r}\cdots\bomega_{n+1-i+s_2, a_2+s_2}\bomega_{n+1-i, a_1} .$$
 Taking $$J=\{n+2-i,\cdots,n\},\ \ \bpi=\bomega_0,\ \ \bpi_1=\bomega_{n+1-i,a_r},\ \ \bpi_2=\bold 1,$$ we see that  Proposition \ref{dsub} gives,
$$\Hom_{\hat\bu_n}(W(\bomega_{n+1-i,\boa}), V(\bomega_2)\otimes V(\bomega_0))\ne 0.$$  
Similarly, taking  $$J=\{1,\cdots ,i-1\},\ \ \bpi=\bomega_0^*,\ \ \bpi_2=\bomega_{i,b_r},\ \ \bpi_1=\bold 1,$$ we get
$$\Hom_{\hat\bu_n}(W(\bomega_{i,\bob}), V(\bomega_0^*)\otimes V(\bomega_1))\ne 0.$$
We claim that in fact we have,
$$\Hom_{\hat\bu_n}(V(\bomega_{n+1-i,\boa}), V(\bomega_2)\otimes V(\bomega_0))\ne 0,\ \ \Hom_{\hat\bu_n}(V(\bomega_{i,\bob}), V(\bomega_0^*)\otimes V(\bomega_1))\ne 0.$$ Assuming the claim we complete the proof of the reverse direction. 
The claim implies that $${\rm{hd}} (V(^*\bomega_2)\otimes V(\bomega_{n+1-i,\boa}))\cong V(\bomega_0)\cong \soc(V(\bomega_1)\otimes V(\bomega_{n+1-i,\boa})),$$
 Using \cite[Theorem 3.12]{KKKO15} this is equivalent to 
$${\soc} (V(\bomega_{n+1-i,\boa})\otimes V(^*\bomega_2))\cong V(\bomega_0)\cong {\rm hd}(V(\bomega_{n+1-i,\boa})\otimes V(\bomega_1)).$$
This means that we have  a non-zero map,$$ V(\bomega_{n+1-i,\boa})\otimes V(\bomega_1)\to V(\bomega_{n+1-i,\boa})\otimes V(^*\bomega_2),$$ and hence a non--zero map,$$f: V(\bomega_1)\otimes V(\bomega_2)\to V(\bomega_{i,\bob})\otimes V(\bomega_{n+1-i,\boa}).$$ 
Since $V(\bomega_1)$ and $V(\bomega_2)$ are not dual representations it follows that $\Im f\ne\mathbb C$. Hence $\Im f$ must contain an $\ell$--highest weight vector with $\ell$--highest weight $\bpi$ satisfying,
 $$\bpi\in\wt_\ell (V(\bomega_1)\otimes V(\bomega_2))\subset \wt_\ell W(\bomega),\ \ \Hom_{\hat\bu_n}(W(\bpi), V(\bomega_{i,\bob})\otimes  V(\bomega_{n+1-i,\boa}))\ne 0.$$
 But now the forward direction of the proposition  implies that $\bpi=\bomega$ which gives the reverse direction of the proposition. 
 \\\\
\noindent It remains to prove the claim. Let $f: W(\bomega_{n+1-i,\boa})\to V(\bomega_2)\otimes V(\bomega_0)$ be a non-zero map and 
 and set $v_0 =f(w_{n+1-i,\boa})$. Then $\hat\bu_n v_0$  is  the image of $f$ and it suffices to prove that   $\hat\bu_n v_0$ is irreducible. 
 For this,   it suffices to show that  $$\bpi\in\wt_\ell^+ W(\bomega_{n+1-i,\boa}),\  \ \Hom_{\hat\bu_n}(W(\bpi), V(\bomega_2)\otimes V(\bomega_0))\ne 0\implies\bpi=\bomega_{n+1-i,\boa}.$$ Setting $s_{r+1}=0$ and $s_1=2i$ we have
$$b_j-s_{j+1}-b_{j-1}-s_j=2i-s_{j+1}-s_j\in S_{i-s_{j+1}, i-s_j, n}. $$ Hence $V(\bomega_s)$, $s=0,1,2$  are  prime snake modules.
By   Lemma \ref{socle} and Proposition \ref{mysnake}  we can write $$\bomega_{n+1-i,\boa}= \bomega_2\bomega(\underline g),\ \ \bpi=\bomega_2\bomega(\underline g'), \ \ \underline g'= (g_2',\cdots, g_r')\in \mathbb P_{\bomegas_0},\ \ \underline g= (g_2,\cdots, g_r)\in \mathbb P_{\bomega_0},$$ where $$\bomega(g_p) = \bomega_{n+1-i+s_p,a_p+s_p}^{-1}\bomega_{n+1-i,a_p}, \ \ 2\leq p\leq r. 
$$
Since $\bomega_{i,\boa}\in\bpi\cal Q^+$, it follows from Proposition \ref{mysnake}(iv)  that 
$$g'(k)\geq g(k), \ \ {\rm{for\ all}}\ \  k\in [0,n+1].$$
Hence if $\underline g'\neq \underline g$ there exists $2\leq p \leq r$ minimal such that $g_p'(n+1-i)>a_p$. 
The definition of $\mathbb P_{n+1-s_p,a_p-i+s_p}$ gives 
$$a_p \leq g_p'(n+1-i)\leq a_p + 2s_p, \ \ \ g_p'(0) = a_p+n+1-i.$$
If $\boc_{g_p'}^-\subset [n+2-i,n]$, since $g_p'(n+1-i)>a_p$ we would have $$g_p'(n+1-i)= g_p'(0)+n+1-i = a_p+2(n+1-i)>a_p+2s_p,$$
and hence a contradiction. It follows that there exists $j\le n+1-i$ with $j\in\boc^-_{g_p'}$ and $g_p'(j)\ge a_p+2$. But this contradicts $\bomega_2\bomega(\underline g')\in\cal P_n^+$. Hence we have proved that $\bpi=\bomega_{n+1-i,\boa}$ as needed. The proof of the second assertion of the claim is identical and we omit the details.\\\\

 \subsection{Imaginary modules in $A_2$}\label{a2}
 In the case of the quantum affine algebra associated to $A_2$, the tensor product of a KR--module and its dual has Jordan--Holder series of length two. Hence the methods of the previous section do not yield an imaginary module. However, the methods do work for a different family of modules and in fact they work in all ranks. \\\\
 The evaluation modules are a particular family of snake modules, which are indexed by elements of $P^+$. In the case of $A_2$ they  correspond to the following elements (and their duals) of $\cal P_2^+$: for $r_1, r_2\in \mathbb N$ and $b_1\in\mathbb Z$, we have 
 \begin{gather}\bpi=\bomega_{1, \bob_1}\bomega_{2,\bob_2},\ \ 
\bob_j=(b_j,b_j+2,\cdots b_j+2r_j-2),\ \ b_1-3=b_2+2r_2-2. \end{gather}  We shall call these elements snakes of evaluation type. Write $${}^*\bpi=\bomega_{1,\boa_1}\bomega_{2,\boa_2},\ \ \boa_1=(b_2-3,\cdots ,b_2+2r_2-5), \ \ \boa_2=(b_1-3,\cdots, b_1+2r_1-5).$$
\begin{prop}
 Suppose that $\bpi$ is a snake of evaluation type in $A_2$  and 
 assume that $r_1, r_2\ge 2$. Then the module $V({}^*\bpi\bpi)$ is imaginary.
\end{prop}
\begin{pf}  Let $M$ be the submodule of $V(\bpi)\otimes V({}^*\bpi)$ generated by the tensor product   $\bold v$ of the highest weight vectors. Consider the quantum affine subalgebra associated to  $J=\{2\}$. It was proved in  \cite{CP91} that  $\hat\bu_{J}\bold v$ is a proper submodule of $V_J(\bomega_{2,\bob_2})\otimes V_J(\bomega_{2,\boa_2})$. It follows that $M$ must be a proper submodule of $V(\bpi)\otimes V({}^*\bpi)$. Next we prove that there exists a non--split short exact sequence \begin{equation}\label{sesm}0\to \mathbb C\to M\to V({}^*\bpi\bpi)\to 0.\end{equation}
For this along of Lemma \ref{dimv},  it is enough to show that
$$\Hom_{\hat{\bu_n}}(W(\tilde\bpi), V(\bpi)\otimes V({}^*\bpi))\ne 0\iff \tilde\bpi=\{\bold 1, {}^*\bpi\bpi\}.$$
Notice that  the reverse direction is obvious.
For the forward direction we use Lemma \ref{socle} to deduce that $\tilde\bpi$ must satisfy,
\begin{gather}\label{equal}\tilde\bpi=\bpi\bomega(\underline g),\ \  {\rm{and}}\ \ ({}^*\tilde\bpi)^{-1}\bpi=\bomega(\underline p),\ \  \underline g\in\mathbb P_{{}^*\bpis},\ \ \ \ \underline p\in\mathbb P_{\bpis}.\end{gather}Writing $$\underline g=(\underline g_1,\underline g_2),\ \ \underline g_1\in\mathbb P_{1,\boa_1}\ \ \underline g_2=(g_1,\cdots, g_{r_1})\in\mathbb P_{2,\boa_2},$$ we see by using Proposition \ref{mysnake} that $$\bomega(g_{r_1})=\bomega_{2,b_1+2r_1-5}\implies\bomega(\underline g_2)=\bomega_{2,\boa_2}.$$ We claim that this further implies that $\bomega(\underline g_1)=\bomega_{1,\boa_1}$. Otherwise we would have an element of the form $\bomega_{1, b_1-4}^{-1} $ or $\bomega_{2,b_1-3}^{-1}$ in a reduced expression for $\bomega(\underline g_1)$. The first would contradict $\tilde\bpi\in\cal P^+$ and the second is impossible by Proposition \ref{mysnake} since $\bomega_{2,b_1-3}$ occurs in $\bomega(\underline g_2)$. This proves the claim and hence we have proved that
$ \tilde\bpi={}^*\bpi\bpi.$\\\\ 
If $\bomega(g_{r_1})\ne \bomega_{2,b_1+2r_1-5}$, then  the condition that $\tilde\bpi\in\cal P^+$ forces $\bomega(g_{r_1})=\bomega_{1, b_1+2r_1-2}^{-1}$ and so we get $$\tilde\bpi=\bomega_{1,\bob_1\setminus\{b_1+2r_1-2\}}\bomega(\underline g_1)\bomega(\underline g_2'),\ \ \underline g_2'=(g_1,\cdots, g_{r-1})\in\mathbb P_{2,\boa_2\setminus\{b_1+2r_1-5\}}.$$
Now we observe that the preceding assertion   implies that $\bomega_{2, b_1+2r_1+1}^{-1}$ does not occur on the left hand side of the second equality in \eqref{equal}. 
Writing $\underline p=(\underline p_2, \underline p_1)$ with $\underline p_s\in \mathbb P_{s,\bob_s}$, $s=1,2$, $\underline p_1=(p_1,\cdots, p_{r_1})$  we see that we must have $\bomega(p_{r_1})=\bomega_{1, b_1+2r_1-2}$.
But this again forces  $\bomega(\underline p_1)=\bomega_{1,\bob_1}$. 
 If $\bomega(\underline p_2)\ne \bomega_{1,\bob_2}$ then an expression of the form $\bomega_{1, b_1-2}$ or $\bomega_{1, b_1-3}^{-1}$ must occur in $\bomega(\underline p_2)$. The first is impossible since such a term does not occur on the left hand side of the second equation in \eqref{equal} while the second would contradict the fact that $\underline p\in\mathbb P_{\bpis}$. Hence we have proved that $\bomega(\underline p_2)= \bomega_{1,\bob_2}$ which now forces $\tilde\bpi=\bold 1$. 
\\\\
 Let $\Phi:V(\bpi)\otimes V({}^*\bpi)\otimes V(\bpi)\otimes V({}^*\bpi)\to V(\bpi)\otimes V({}^*\bpi) $ be defined as in Section \ref{phi}. Again Lemma \ref{dimv} shows that $\Phi(M\otimes M)$ is a  non--zero submodule of $V$. By Proposition \ref{lowestwt} we have $$\Phi(M_{(\bpis\bpis^*)^{-1}}\otimes M_{{}^*\bpis\bpis})\ne 0.$$ A standard argument (as in the previous examples) using Proposition \ref{mysnake} shows that  $\dim (V(\bpi)\otimes V({}^*\bpi))_{{}^*\bpis(\bpis^*)^{-1}}=1$  and so we get,
 $$\Phi(M_{(\bpis\bpis^*)^{-1}}\otimes M_{{}^*\bpis\bpis})=(V(\bpi)\otimes V({}^*\bpi))_{{}^*\bpis(\bpis^*)^{-1}}.$$
Let $\bar\bpi\in\cal P^+$ be the head $V(\bpi)\otimes V({}^*\bpi)$ and note that $\bar\bpi\ne \bold 1,{}^*\bpi\bpi, ({}^*\bpi\bpi)^2$. 
Lemma \ref{lowestwt} now gives that the image of $(V(\bpi)\otimes V({}^*\bpi))_{{}^*\bpis(\bpis^*)^{-1}}$ in  $V(\bar\bpi)$ is non--zero. The discussion so far shows that the composite map:
$$M\otimes M\to V(\bpi)\otimes V({}^*\bpi)\to V(\bar\bpi)$$ is surjective.  Using the short exact sequence in \eqref{sesm} we see that the composite map obviously factors through to a map from $V({}^*\bpi\bpi)\otimes V({}^*\bpi\bpi)  $ and the proof of the proposition is complete.
 \end{pf} 
  \subsection{$D_4$}\label{d4}  We give the first family of  examples   of \ imaginary modules in type $D_4$ which does not arise from an embedding of $A_3$. Assume that $2$ is the trivalent node in $D_4$. We shall prove that if $$\bomega=\bomega_{1,2r}\bomega_{1,2r-2}\bomega_{2,2r-5}\cdots\bomega_{2,5}\bomega_{1,2}\bomega_{1,0}$$ then the module $V(\bomega)$ is imaginary $r\geq 5$.  The proof follows the same strategy as in the case of $A_n$.
  \begin{rem}Note that regarded as module for the quantum affine $A_3$ subalgebras given by $\{j,2,k\}$ with $j\ne k\in\{1,3,4\}$ the module $V(\bomega)$  is an example of a snake module \cite{MY12} and hence  by \cite{DLL19} real for the subalgebra. \hfill\qedsymbol
  \end{rem}
  It  follows from \cite{CP95} that  $V(\bomega_{j,a})^*\cong V(\bomega_{j,a+6})$.
We take $$\bob=(6,8,\cdots ,2r), \ \  \boa= (0,2,\cdots, 2r-6),$$ $$ V=V(\bomega_{1,\bob})\otimes V(\bomega_{1,\boa}). $$ 
 We  have  a map of $\hat \bu(D_4)$--modules
 $\mathbb C\hookrightarrow V$ and a projection $<\ ,\ >:V(\bomega_{1,\boa})\otimes V(\bomega_{1,\bob})\to\mathbb C$ and we let $\Phi: V\otimes V\to V$ be the composition  $\id\circ<\ ,\ >\circ\id$. \\\\ It is helpful to recall that  the elements of $\wt_\ell V(\bomega_{1,0})$ are the following: \begin{gather*}\bomega_{1,0},\ \bomega_{2,1}\bomega_{1,2}^{-1},\ \bomega_{2,3}^{-1}\bomega_{3,2}\bomega_{4,2},\ \bomega_{3,4}^{-1}\bomega_{4,2},\  \bomega_{4,4}^{-1}\bomega_{3,2},\\ \bomega_{2,3}\bomega_{3,4}^{-1}\bomega_{4,2}^{-1},\ \ \bomega_{1,4}\bomega_{2,5}^{-1},\ \ \bomega_{1,6}^{-1}. \end{gather*}
 Notice that $$\bomega=\bomega_{i,\bob}\bomega',\ \ \bomega'=\bomega_{i,\boa}\balpha_{1,2r-5}^{-1}\cdots\balpha_{1,5}^{-1}\in\wt_\ell V(\bomega_{i,\boa}).$$
 \subsubsection{} The following result will be used.
 \begin{lem}\label{cornerd}
 Suppose that $\bpi\in\wt_\ell V(\bomega_{i,\boa})$ is such that  $\bpi=\bomega_{1,c_1}^{-1}\cdots\bomega_{1,c_r}^{-1}\bpi'$ for some $c_1,\cdots, c_k$ and  $\bpi'\in\cal P^+$. There exists $1\le s\le p\le r+1$ such that  $$\bpi=\bomega_{1,a_r+6}^{-1}\cdots\bomega_{1,a_{p}+6}^{-1}(\bomega_{1, a_{p-1}+2}^{-1}\bomega_{2,a_{p-1}+1})\cdots(\bomega_{1,a_{s}+2}^{-1}\bomega_{2,a_{s}+1})\bomega_{1,a_{s-1}}\cdots\bomega_{1,a_1}.$$
 \end{lem}
 \begin{pf} Since $\bpi\in\wt_\ell V(\bomega_{i,\boa})$ and $\bpi\ne \bomega_{i,\boa}$ it follows that there exists $1\le s\le r$ minimal such that $$\bpi\le \bomega_{i,\boa}\balpha_1^{-1},\ \ \balpha_1=\balpha_{1,a_r+1}\cdots\balpha_{1,a_s+1}.$$ 
 Hence either the Lemma is proved or the inequality is strict. In the latter case, there must exist $p\ge s$ minimal  with $$\bpi\le  \bomega_{i,\boa}(\balpha_1\balpha_2)^{-1},\ \ \balpha_2=\balpha_{2,a_r+2}\cdots\balpha_{2, a_{p}+2}.$$ The assumptions on $\bpi$ now means that this inequality must be strict. Now working with $\balpha_j$ with $j=3,4 $ we see in a similar fashion that  $$\bpi\le \bomega_{i,\boa}(\balpha_1\balpha_2\balpha_3\balpha_4)^{-1},\ \ \balpha_j=\balpha_{j,a_r+3}\cdots\balpha_{j,a_{p}+3},\ \ j=3,4.$$ Again the assumption on $\bpi$ mean that the inequality must be strict and so we get $$\bpi\le \bomega_{i,\boa}(\balpha_1\balpha_2\balpha_3\balpha_4\balpha_2')^{-1},\ \ \balpha_2'=\balpha_{2,a_r+4}\cdots\balpha_{2,a_{p}+4}.$$ A final application of the property of $\bpi$ gives$$\bpi\le \bomega_{i,\boa}(\balpha_1\balpha_2\balpha_3\balpha_4\balpha_2'\balpha_1')^{-1},\ \ \balpha_1'=\balpha_{1,a_r+5}\cdots\balpha_{1,a_{p}+5},$$ or equivalently
 $$\bpi\le \bomega_{1,a_r+6}^{-1}\cdots\bomega_{1,a_{p}+6}^{-1}(\bomega_{1, a_{p-1}+2}^{-1}\bomega_{2,a_{p-1}+1})\cdots(\bomega_{1,a_s+2}^{-1}\bomega_{2,a_s+1})\bomega_{1,a_{s-1}}\cdots\bomega_{1,a_1}.$$ The minimality of the choice of $p,s$ now proves this must be an equality as needed.
 \end{pf}
 
 \subsubsection{} \begin{lem}
Let $\bvarpi\in\wt_\ell^+W(\bomega)$. Then $$\Hom_{\hat\bu}(W(\bvarpi), V)\ne 0\iff \bvarpi\in\{\bold 1,\bomega\}.$$
 \end{lem}
 \begin{pf} 
 The converse is clear by using the result for quantum affine $\lie{sl}_2$. 
 For the forward direction we use Lemma \ref{socle} to write 
 $\bvarpi=\bomega_{1,\bob}\bpi$ with $\bpi\in\wt_\ell V(\bomega_{1,\boa})$. Since $\bvarpi\in\cal P_n^+$ we see that $\bpi$ must satisfy the conditions of Lemma \ref{cornerd} and let $s,p$ be as in that Lemma. \\\\
 Since $\bvarpi\le\bomega$ we are forced to take $s\in\{1,2,3\}$.
 If $s=1,2$ and if $p-1\ge s$ then the term $\bomega_{1,a_2+2}^{-1}$ occurs in a reduced expression for  $\bvarpi\in\cal P^+$ which is absurd. Hence 
 we must have $p-1<s$ which forces $p=s$. This means that $\bvarpi\in\{\bold 1, \bomega_{1,b_1}\bomega_{1,a_1}\}$. Using the second assertion of Lemma \ref{socle} we see that $$\bvarpi=\bomega_{1,b_1}\bomega_{1,a_1}\implies \bomega_{1,b_1+6}^{-1}\bomega_{1,b_1}^{-1}\bomega_{1,\bob}\in\wt_\ell V(\bomega_{1,\bob}).$$  
 If $r\geq 6$ we get a contradiction since $\bomega_{1,b_1+6}^{-1}\bomega_{1,b_1}^{-1}\bomega_{i,\bob}\in \cal P^+$ and $\wt_\ell^+ V(\bomega_{1,\bob})=\{\bomega_{1,\bob}\}$, by \cite{Nak04} (see also \cite[Theorem 2.12]{He06}). If $r=5$ this means that 
 $$\bomega_{1,b_2}\bomega_{1,b_3}\bomega_{1,b_1+6}^{-1}\in \wt_\ell V(\bomega_{1,\bob}),$$
  It follows from \cite[Lemma 5.8]{He06} that %$$\bomega_{1,b_2}\bomega_{1,b_3}\bomega_{1,b_1+6}^{-1} \in  \bomega_{1,\bob}(\balpha_{1,b_1+1}\balpha_{1,b_2+1}\balpha_{1,b_3+1})^{-1} (\cal Q^+)^{-1} \implies 
 $$\bomega_{1,b_2}\bomega_{1,b_3}\bomega_{1,b_1+6}^{-1} \in (\bomega_{1,b_1+2}\bomega_{1,b_2+2}\bomega_{1,b_3+2})^{-1} \bomega_{2,b_1+1}\bomega_{2,b_2+1}\bomega_{2,b_3+1} (\cal Q^+)^{-1},$$
 which is absurd since $\bomega_{1,b_2}$ can never occur in the right hand side. 
% {\color{red} Write $$\bomega_{1,b_1+6}^{-1}\bomega_{1,b_1}^{-1}\bomega_{i,\bob}=\bomega_1\cdots\bomega_r,\ \ \bomega_s\in\wt_\ell V(\bomega_{1,b_s}).$$  Since $r\ge 5$ this means that $\bomega_{i,b_r}$ definitely occurs in  $\bomega_{1,b_1+6}^{-1}\bomega_{1,b_1}^{-1}\bomega_{i,\bob}$. This forces $\bomega_s=\bomega_{1,b_r}$ and hence by \cite{He05} that $\bomega_s=\bomega_{1,b_s}$ for all $s$ which is a contradiction.}
 \\\\
 Suppose next that $s=3$ and let $3\leq p\leq r+1$. If $p=r+1$ then $\bvarpi=\bomega$ and there is nothing to prove.
 If $p\leq r$ then,
 $$\bvarpi = \bomega_{1,a_1}\bomega_{1,a_2}\bomega_{2,a_3+1}\cdots \bomega_{2,a_{p-1}+1}\bomega_{1,b_{p-2}}\bomega_{1,b_{p-1}}.$$
 Using Lemma \ref{socle} we have
 %we have $$\left(\bomega_{1,b_1}\bomega_{1,b_2}\bomega_{2,b_3+1}\cdots \bomega_{2,b_{p-1}+1}\bomega_{1,b_{p-2}+6}\bomega_{1,b_{p-1}+6}\right)^{-1}\bomega_{1,\bob}\in\wt_\ell V(\bomega_{i,\bob}).$$
% If $p\le r-3$ or if $p=r$ then $\bomega_{i,b_r}$ occurs in the reduced expression for $\bvarpi$ and  the argument given above works. \\\\
 $$(\bomega_{1,b_3}\bomega_{2,b_3+1}^{-1})\cdots (\bomega_{1,b_{p-1}}\bomega_{2,b_{p-1}+1}^{-1}) (\bomega_{1,b_{p-2}+6}\bomega_{{1,b_{p-1}+6}})^{-1} \bomega_{1,b_p}\cdots \bomega_{1,b_r}\in \wt_\ell V(\bomega_{1,\bob}).$$
 In particular we can write the above element as 
 $$\bomega_{1}\cdots \bomega_{p-1}\bomega' \in \wt_\ell V(\bomega_{1,b_1})\otimes \cdots \otimes  V(\bomega_{1,b_{p-1}})\otimes  V(\bomega_{1,\bob'}), \ \ \bomega'\in \wt_\ell V(\bomega_{1,\bob'}),$$
 where  $\bob' = (b_p,\cdots,b_r)$. 
 Comparing the possible weights we are forced to have
 $$\bomega_s = \bomega_{1,b_{s+2}}\bomega_{2,b_{s+2}+1}^{-1}, \ \ \ 1\leq s\leq p-3,$$
 and hence 
 $$(\bomega_{1,b_{p-2}+6}\bomega_{{1,b_{p-1}+6}})^{-1} \bomega_{1,b_p}\cdots \bomega_{1,b_r}\in \wt_\ell V(\bomega_{1,\bob'}).$$
 If $p\leq r-2$ the corresponding $\ell$--weight is domintant and again we have a contradiction since  $\wt_\ell^+V(\bomega_{1,\bob'})=\{\bomega_{1,\bob'}\}$. If $p=r-1$ it means that
 $$\bomega_{1,b_p}\cdots \bomega_{1,b_{r-1}}\bomega_{1,b_{r}+2}^{-1}\in \wt_\ell V(\bomega_{1,\bob'})\implies \bomega_{1,b_{r}+2}^{-1}\in \wt_\ell V(\bomega_{1,b_r}),$$
 Similarly, if $p=r$ we have
 $$\bomega_{1,b_p}\cdots \bomega_{1,b_{r}}\bomega_{1,b_{r}+2}^{-1}\bomega_{1,b_{r}+4}^{-1}\in \wt_\ell V(\bomega_{1,\bob'}).$$
 In either case we have a contradiction. 
 \end{pf}

 \subsubsection{}  The proof that $V(\bomega)$ is imaginary is now completed as in the case of $A_n$. Let $f: W(\bomega)\to V$ be a non--zero map of $\hat\bu$--modules and let $v=f(v_\bomegas)$ and set $M(\bomega)=\hat\bu v$. The preceding Lemma and \cite{KKKO15} show that we have a non--split short exact sequence $$0\to \mathbb C\to M(\bomega)\to V(\bomega)\to 0.$$
  Suppose that $\bomega_1=\bomega_2^{-1}$ with $\bomega_1\in \wt_\ell V(\bomega_{1,\bob})$
 and $\bomega_2\in \wt_\ell V(\bomega_{1,\boa})$. Using \cite{He06} again we see that if $\bomega_1\ne \bomega_{1,\bob}$ then $\bomega_1$ must involve a term $\bomega_{j,c}$ with $c>b_r$. Since $\bomega_{j,c}^{-1}$ cannot occur in an element of $\wt_\ell V(\bomega_{1,\boa})$ we see that $\bomega_1=\bomega_{1,\bob}$ and so $\dim V_{\bold 1}=1$.  Similarly if $\bomega_1\bomega_2=\bold 1$ with $\bomega_1\in \wt_\ell V(\bomega_{1,\bob})^{\otimes 2}$
 and $\bomega_2\in \wt_\ell V(\bomega_{1,\boa})^{\otimes 2}$ one proves that $\bomega_1=\bomega_{1,\bob}^2$ and hence $\dim (V\otimes V)_{\bold 1}=1$. It follows that the restriction of $\Phi$ to $M(\bomega)\otimes M(\bomega)$ is non--zero.
 \\\\
 Next we note that $\bomega^2\notin \wt_\ell V$. This proves that $V(\bomega^2)$ cannot be in the head of $\Phi(M(\bomega)\otimes M(\bomega)).$ 
 \\\\
 Finally, we  prove that $(\bomega^*)^{-1}\bomega\notin\wt_\ell W(\bomega)$.  Suppose that $(\bomega^*)^{-1}\bomega=\bomega_1\cdots\bomega_{2r+2}$ where $$\bomega_s\in\wt_\ell W(\bomega_{1, 2s}),\  s=r-1, r, \ \ \bomega_s\in W(\bomega_{1,2s}),\  s=1,2,\ \ \bomega_s\in\wt_\ell W(\bomega_{2, 2s-1})\ \ {\rm{otherwise}}. $$  Then we must have $\bomega_{s}=\bomega_{1, 2s}$ if $s=1,2$.  But then $\bomega_{1,6}^{-1}$ does not occur as part of an $\ell$--weight of the other modules and we are done.  Hence Proposition \ref{lowestwt} shows that $V(\bomega)$ cannot occur in the head of $\Phi(M(\bomega)\otimes M(\bomega)).$  Since $\bomega$ is not self dual it is now clear that $\Phi(M(\bomega)\otimes M(\bomega))$ has an irreducible module $V(\bpi)$ in its head with $\bpi\ne\bold 1,\bomega,\bomega^2$ and this map factors through to give a non-zero map $V(\bomega)\otimes V(\bomega)\to V(\bpi)$ proving that $V(\bomega)$ is imaginary.

 \end{document}